\documentclass[aip,jcp,preprint,amsmath,amssymb]{revtex4-1}
\usepackage{bm,bbm}
\usepackage{graphicx}
\usepackage{url}

\usepackage{color}

\def\Ee{\mathcal{E}}

\def\Gg{\mathcal{G}}

\def\Ll{\mathcal{L}}

\def\Vv{\mathcal{V}}

\def\E{\mathbb{E}}

\def\N{\mathbb{N}}

\def\R{\mathbb{R}}

\def\EXPECT{{\mathbb{E}}}

\def\RELENT#1#2{\Rr\left(#1|#2\right)}

\def\COMMA{\,,}             
\def\PERIOD{\,.}            
\def\SEP{{\,|\,}}           

\def\VIZ#1{(\ref{#1})}      

\def\BIGO{\mathcal{O}}

\def\PROOF{\noindent{\sc Proof:}\ }

\newtheorem{thm}{Theorem}
\newtheorem{lem}[thm]{Lemma}

\newtheorem{cor}{Corollary}

\def\0o{\mathbb{O}}

\newcommand\ontop[2]{\genfrac{}{}{0pt}{}{#1}{#2}}

\def\bx{z}

\def\Id{\mathbf{I}}
\def\BFW{\mathbf{W}}
\def\BFT{\mathbf{T}}
\def\BFD{\mathbf{D}}
\def\BFJ{\mathbf{J}}
\def\BFG{\mathbf{G}}

\def\BARIT#1{{\bar {#1}}}

\def\BARM{\BARIT{\mu}}

\def\BARV{\BARIT{U}}
\def\IP{U}
\def\PMF{\BARIT{U}^{\mathrm{PMF}}}
\def\FPMF{F^{\mathrm{PMF}}}

\def\BARMV{\BARIT{\mu}_{\BARIT U}}

\def\bm{\BARIT{m}}

\def\COP{\mathbf{\Pi}}

\def\BIGO{\mathcal{O}}
\def\VIZ#1{(\ref{#1})}

\def\SEP{\,|\,}

\def\RELENT#1#2{\mathcal{R}\left({#1}\SEP{#2}\right)}

\def\RELENT#1#2{\mathcal{R}\left({#1}\SEP{#2}\right)}

\def\SEP{{\,||\,}}           

\def\BARIT#1{{\bar {#1}}}

\def\BARM{\BARIT{\mu}}

\def\COP{\mathbf{\xi}} 
\def\COPL{\mathbf{T}}                

\def\EXPECT{{\mathbb{E}}}

\begin{document}
\title{The geometry of generalized force matching in coarse-graining and related information metrics}
\author{Evangelia Kalligiannaki}
\email{ekalligian@tem.uoc.gr}
\affiliation{Department of Mathematics and Applied Mathematics, University of Crete, 70013 Heraklion, Greece}
\author{Vagelis Harmandaris}
\email{harman@uoc.gr}
\affiliation{Department of Mathematics and Applied Mathematics, University of Crete, 70013 Heraklion, Greece}
\affiliation{  IACM/FORTH GR-71110 Heraklion, Greece}
\author{Markos A. Katsoulakis}
\email{markos@math.umass.edu}
\affiliation{Department of Mathematics and Statistics, University
                of Massachusetts,  Amherst, MA 01003, USA}
\author{Petr Plech\'a\v{c}}
\email{plechac@math.udel.edu}
\affiliation{Department of Mathematical Sciences,
           University of Delaware,
           Newark, DE 19716, USA }
\date{\today}

\begin{abstract}
 Using the probabilistic language of conditional expectations we reformulate the force matching method for
coarse-graining of molecular systems as a projection on spaces of coarse observables.
A practical outcome of this probabilistic description  is the link  of the
force matching method   with thermodynamic integration.
This connection provides a way to systematically construct  a
local mean force in order to optimally approximate the potential of mean force  through  force matching.
We introduce a  generalized force matching condition  for the local mean force
in the sense that allows the approximation of the potential of mean force under both linear and
non-linear coarse graining    mappings (e.g., reaction coordinates, end-to-end length of chains).
Furthermore, we study the equivalence of force matching with relative entropy minimization which we
derive for general non-linear coarse graining maps. We present in detail the generalized  force matching
condition  through applications to specific examples  in  molecular systems.
 
\end{abstract}

\maketitle

\section{Introduction}
Complex molecular systems are materials of amazing diversity ranging from polymers to colloids, hybrid nanocomposites, biomolecular systems, etc, which are directly related with an enormous range of possible applications in nano-, bio-technology, food science, drug industry, cosmetics etc. Due to the above reasons molecular simulations of complex systems is a very intense research area.\cite{FrenkelSmitBook} A main challenge in this field is to predict structure-properties relations of such materials and to provide a direct quantitative link between chemical structure at the molecular level and measurable structural and dynamical quantities over a broad range of  length and time scales.

On the microscopic (atomistic) level, detailed all-atom molecular dynamics (MD), or Monte Carlo (MC) simulations allow direct quantitative predictions of the properties of   molecular systems over a range of length and time scales. \cite{AllenTildesleyBook,FrenkelSmitBook,Harmandaris2003a} However, due to the broad spectrum of characteristic lengths and times involved in complex molecular systems it is not feasible to apply them to large realistic systems or molecules of complex structure, such as multi-component biomaterials, polymers of high molecular weight, colloids etc. 
On the mesoscopic level, coarse-grained (CG) models have proven to be very efficient means in order to increase the length and time scales accessible by simulations.\cite{FrenkelSmitBook,IzVoth2005a,tsop1,MulPlat2002,Shell2008,briels,Harmandaris2003a,Harmandaris2006a,Harmandaris2009a,
Harmandaris2009b,Johnston2013,IzVoth2005a,Voth2008a,Voth2010,Noid2011,Noid2013,Shell2009,Zabaras2013,DiffusionMapKevrekidis2008,Soper1996,LyubLaa2004}

CG (particle) models can be roughly categorized, based on the way they are developed, into two groups: (a) Ad hoc or phenomenological CG models, such as simple bead spring or lattice ones, which are primarily used to study generic behavior (e.g. scaling properties) of complex systems but lack a link to specific systems.\cite{FrenkelSmitBook} The interactions between the CG groups in these models are described through semi-empirical functional forms obtained through previous knowledge and with a lot of physical intuition. 
(b) Systematic CG models, which are usually developed by lumping groups of atoms into groups, i.e. "superatoms", and deriving the effective CG interaction potentials directly from more detailed (microscopic) simulations. Such models are capable of predicting \textit{quantitatively} the properties of \textit{specific systems} and have been applied with great success to a very broad range of molecular systems (see for example refs. \cite{IzVoth2005a,tsop1,MulPlat2002,Shell2008,briels,Harmandaris2009a,Harmandaris2009b,Johnston2013} and references within).
 
A main challenge in the later family of CG models is to develop rigorous atomistic to CG methodologies that allow, as accurate as possible,  the estimation of the CG effective interaction. With such approaches the hierarchical combination of atomistic and CG models could be in order to study a very broad range of length and time scales of specific molecular complex systems without adjustable parameters, and by that become truly predictive. \cite{Voth2008a,Harmandaris2009a,Voth2010}

Let us assume a specific molecular system. The overall procedure of systematic coarse-grained modeling for this system, based on detailed microscopic data, is shortly described through the following stages:
(a) Execution of microscopic (e.g. ab-initio or atomistic) simulations on small model systems, i.e. usually a relatively small number of molecules with a rather low molecular weight is considered, (b) Choose of the CG map (transformation from the atomistic to the CG description), (c) Development of the CG effective interaction potential (force field), (d) Execution of the CG (e.g. MD, Langevin dynamics, LD, or MC) simulations and (e) Re-introduction of the missing atomistic degrees of freedom in the CG structures, in case the properties under study require atomistic detail.
From all above stages the development of the CG force field is the most challenging one. Indeed, an accurate estimation of the way CG "superatoms" interact to each other is a \textit{conditio sine qua non} in order to understand the behavior and to (quantitatively) predict the properties of the specific complex molecular system under study. 

Note that from a mathematical point of view coarse-graining is a sub-field of the dimensionality reduction.\cite{pavliotis2008} Indeed, there are several statistical methods for the reduction of the degrees of freedom under consideration, in a deterministic or stochastic model, such as principal component analysis and diffusion maps.\cite{DiffusionMapKevrekidis2008} Here we focus our discussion on CG methods based  on statistical mechanics, which are used extensively the last two-three decades in the theoretical modeling of molecular systems across a very broad range of disciplines, from physics to chemistry and biology as well as in engineering sciences.     

There exists a variety of methods that construct a  reduced  model that approximates effective properties of complex systems based on statistical mechanics. These methods usually consider the optimization of proposed parametric models using different  minimization principles, that is considering a pre-selected set of observables $\{\phi_i, \ i=1,\dots,k\}$ and then minimizing over a parameter set $\Theta$,
\begin{equation*}
\min_{\theta\in \Theta} \sum_{i=1}^k \|\E_{\mu}[\phi_i] - \E_{\mu^{\theta}} [\phi_i]\|^2 \COMMA
\end{equation*}
where $\mu(x), \ \mu^{\theta}(x)$ are the atomistic and proposed Gibbs measures respectively.
Different methods consider different sets of observables. For example:

\noindent
(a) In structural based methods the observable  is the {\it pair radial distribution function} $g(r)$, related to the two-body potential of mean force (see section~\ref{IBIM}), for the intermolecular interaction potential, and distribution functions of bonded degrees of freedom (e.g. bonds, angles, dihedrals) for CG systems with intramolecular interaction potential.\cite{Soper1996,LyubLaa2004,tsop1,MulPlat2002,Harmandaris2006a} 
  
 \noindent 
(b) Force matching (FM) methods \cite{IzVoth2005,IzVoth2005a,Noid2011} consider as observable function     the force $f_{j}(x)=-\nabla_{x_j} \IP(x),\ j=1,\dots,N$, for an $N$-particle system with interaction potential $\IP(x),\  x\in \R^{3N}$.

\noindent
(c) The relative entropy (RE)\cite{Shell2008,Shell2009,Zabaras2013}  method employs the minimization  of the  relative entropy  pseudo-distance 
$$
   \RELENT{\mu}{ \mu^{\theta}} = \int_{\R^{3N}} \log \frac{d\mu(x)}{d\mu^{\theta}(x)} d\mu(x)\COMMA
 $$
These  methods, in principle, are employed to  approximate a many body potential   describing the {\em equilibrium} distribution of  CG particles observed in simulations of atomically detailed models. 
The many body potential is  defined through the renormalization group map \cite{Goldenfeld1992} that is equivalent to the potential  of mean force (PMF) \cite{mcquarrie2000statistical} in case  the former is differentiable.  
The  force-matching (or multi scale coarse graining  (MSCG) )  and the relative  entropy are minimization methods that   construct a best fit of  a proposed coarse graining    potential for systems in {\em equilibrium}. The force-matching method determines a CG potential  from  {\it  atomistic force  information} through a least-square minimization principle, to variationally  project the force   corresponding to the potential of mean force onto a force   that is defined by the form of the approximate potential. The relative entropy approach  obtains optimal CG potential parameters by minimizing the relative entropy between the atomistic and the CG {\it Gibbs measures} sampled by the atomistic.    
A  brief review and categorization of parametrization  methods in equilibrium is given in ref. \cite{Noid2013}

Besides all the  above, a classical method for calculating free energy differences using arbitrary reaction coordinates is thermodynamic integration (TI) theory.\cite{Sprik1998, denOtter, denOtterBriels, Ciccotti} Thermodynamic integration is based on writing free energy differences as the integral of  free energy derivative  and  thus computing  the derivatives (mean force) instead of directly the free energy.
     
The purpose of this work is: 
(a)   To reformulate in the {\it probabilistic language of conditional expectations} the force matching method. 
In turn, the conditional expectation formulation allows us: (b) To reveal the connection of force matching  with thermodynamic integration that provides  a way to {\it construct a local mean force} in order to best approximate the potential of mean force when applying the  force matching method. 
(c) To present in a   probabilistic formalism  the    equivalence of  relative entropy  and force matching methods  which we  derive for general {\it nonlinear coarse graining maps}. 
  We  furthermore discuss structure-based (SB) CG methods 
   thus presenting  a complete picture of the known many body potential estimation methods for systems at equilibrium and their relation.

 Furthermore, the probabilistic formalism gives a geometric representation of the force matching method, i.e.   recast the force matching as a projection procedure onto the space of coarse obsvervables (we refer specifically to Figure~\ref{scheme2}).
The novelty and advantages of our approach is that it allows us to define a generalized force matching  minimization problem 
$$\min_G \E_\mu[\|h - G(\xi)\|^2] $$
applicable  for linear and \textit{nonlinear CG maps}  $\xi:\R^{3N}\to \R^{m}$.
 The  force matching condition introduced  
  \begin{equation}\label{hx}
h(x) =  \BFJ{\COP}^{-1}(x) \BFD\COP(x)  f(x) +\frac1\beta \nabla_{x} \cdot  \BFJ{\COP}^{-1}(x)\BFD\COP(x)\COMMA
 \end{equation}
 where  $ \BFJ{\COP}(x) = \BFD\COP(x)\BFD\COP^t(x) $,$\BFD\COP\in \R^{m\times 3N}$ $(\BFD\COP)_{ij}(x) =  \nabla_{x_j} \COP_i(x) , \ i=1,\dots,m, j=1,\dots,N$
 ensures the best approximation of the  PMF. As in thermodynamic integration $h(x)$ is called the {\it local mean force}. A more general result  is available in Theorem~\ref{ghPMF}.
This elucidates the direct connection of the above discussed particle CG methodologies with the standard thermodynamic integration approaches.

The current work is directly related to previous works that concern linear CG mapping schemes.\cite{Noid2011, Voth2008a} Here we recast and extend these works in a probabilistic formalism in order to  present  and compare the relative entropy and force matching methods that allows us to generalize the methodology to nonlinear coarse-graining maps.
 In the case of a linear CG map the local mean force  is
$$ h(x)=  \BFJ{\COP}^{-1}(x) \BFD\COP(x)  f(x) $$
which reduces to the result in ref. \cite{ Voth2008a}, see the examples in Section~\ref{LinearEx}. Notice that for linear CG maps   the last term in relation \VIZ{hx} vanishes. 
The proposed formula for the local mean force extends the  works\cite{Noid2011, Voth2008a} in two aspects: (a) to any non-linear CG map  and (b) the  existence of a family of appropriate local mean forces $h(x)$.



Finally, we should note that in the above discussion we focus on molecular systems at equilibrium. The study of non-equilibrium systems is an even more challenging area related to various phenomena such as the response of the molecular systems on external stimuli (e.g. rheological properties or mechanical behavior of composites systems). The development of CG force fields for systems under non-equilibrium conditions, based on the information theory and  path-space tools of relative entropy is the subject of ref.\cite{HKKP}.

The structure of this work is as follows. In Section~\ref{AtCG} we introduce the atomistic molecular system and its coarse graining through the definition of the CG map and the $n$-body potential of mean force. 
The probabilistic formulation   
  of the force matching method  and  the  best approximation of the PMF are presented in  Section~\ref{CondFM}, while the force matching  condition for approximating the PMF for any, linear or non-linear CG map, are given in Section~\ref{GenFM}. In Sections~\ref{LinearEx}  and~\ref{ReactionC}  we calculate the analytic form of the local mean force for examples of  linear and non-linear CG maps  in molecular systems. 
A second result of the current work is presented in Section~\ref{REandFM} where we prove that the relative entropy minimization and the force matching methods are equivalent, producing the same approximation to PMF up to a constant. Furthermore, for completeness   we present the structure based methods in Section~\ref{IBIM}. We close with Section~\ref{Conclude} summarizing and discussing the results of this work.
\section{Atomistic and coarse-grained systems}\label{AtCG}
Assume the prototypical problem of $N$ (classical) molecules in a box of volume $V$ at temperature $T$.
Let $x=(x_1,\dots,x_N) \in \R^{3N}$ describe the position vectors of the $N$ particles in the atomistic (microscopic)  description, with 
  potential  energy $\IP(x)$.  The probability of a state $x$ at temperature $T$ is given by the Gibbs canonical measure
$$ \mu(dx) =Z^{-1}\exp\{ -\beta \IP(x)\}dx\COMMA$$ 
where $Z= \int_{\R^{3N}}e^{-\beta \IP(x)} dx$ is  the partition function, $\beta=\frac{1}{k_B T}$ and  $k_B$ the Boltzmann constant.  We denote  $f(x)$ the force corresponding to the potential $\IP(x)$  that is $ f: \R^{3N}\to \R^{3N} \nonumber $
\begin{eqnarray}
 f_j(x) = -\nabla_{x_j} \IP(x), \ \ j=1,\dots,N\PERIOD
\end{eqnarray}
 
For such a system  the  $n$-body, $n<N$, potential of mean force (PMF) $ \PMF\!(x_1,\dots,x_n)\!$ \cite{Kirkwood1935, mcquarrie2000statistical} is defined through
  $$  \PMF (x_1,\dots,x_n) = -\frac{1}{\beta} \log g^{(n)}(x_1,\dots,x_n) \COMMA $$ 
where $g^{(n)}(x_1,\dots,x_n)$ is the $n$-body distribution function
 \begin{equation*}
 	g^{(n)}(x_1,\dots,x_n) = \frac{N!}{(N-n)! \rho^n}\int_{\R^{N-n}} \mu(x) dx_{n+1}\dots dx_N\COMMA
 \end{equation*}
 and $\rho=\frac{N}{V}$ is the number density.  

Coarse-graining is considered as the application of a  mapping (CG mapping)   $\COP:\R^{3N} \to \R^{3M}$
\begin{eqnarray}
   x  \mapsto \COP (x)\in \R^{3M}\nonumber
\end{eqnarray}    
 on the microscopic state space,  determining the $M(<N)$ CG particles as a function  of the atomic configuration $x$. 
  We denote by $  \bx = (\bx_1,\dots,\bx_M)$ any point in the CG configuration space $\R^{3M}$, and use the $\BARIT{\cdot }$ notation for quantities on the CG space. We call 'particles' the elements of  the microscopic space with positions $ x_j\in \R^3, j=1,\dots,N $ and 'CG particles' the elements of the coarse space  with positions $ \bx_i\in \R^3,\ i=1,\dots,M$.
We should note that a  CG mapping $\xi(x)$ does not necessarily   maps to  three-dimensional CG particles, it can be considered in the more general form $\xi:\R^{3N}\to \R^m$, for any $m\in \N$ with $m<3N$.  This is the case, for example, when considering some reaction coordinates, like the bending angle and the end-to-end distance (see examples in Section~\ref{ReactionC}).
 
The conditional  Helmholtz free energy $A(\bx)$ related to the CG mapping $\COP(x)$, defined by the renormalization group map\cite{Goldenfeld1992}, is 
 based on the property that for any observable $\phi: \R^{3N}\to \R$ of the form $\phi (x) = \phi(\COP (x))$ it holds 
 \begin{eqnarray*}
 \E_\mu[\phi] &=&\int_{\R^{3N}} \phi(\COP(x)) \mu(x) dx = \int_{\R^{3M}} \int_{\{x:\ \COP( x) = \bx\}}  \phi(\bx) \mu(x) dx d\bx \\
 &=& \int_{\R^{3M}}\phi(\bx)  \BARM(\bx) d\bx = \E_{\BARM}[\E_{\mu}[\phi|\COP]]\COMMA
 \end{eqnarray*}
where $\E_\mu[\phi]  $ denotes the expectation of $\phi(x)$ with respect to the probability measure $\mu(dx)$, 
  \begin{equation}\label{exactCGm}
  \BARM(\bx) =  \int_{\Omega(\bx)}\mu(x)dx, 
  \ \ \Omega(\bx) =\{x\in \R^{3N}: \ \COP(x) = \bx\}  \COMMA
   \end{equation}
 and $\E_\mu\left[\phi|\COP\right]$ is defined by 
  \begin{equation}\label{CE}
  \E_\mu\left[\phi|\bx\right]:= \E_\mu\left[\phi|\COP=\bx\right] =  \int_{\Omega(\bx) }  \phi(x) \mu(x)dx, \ \text{ for all } \bx\in \R^{3M} \  \COMMA
   \end{equation}
   being  the  conditional expectation of  the observable quantity  $\phi(x) $ that represents the expectation of $\phi(x)$ with respect to the Gibbs measure $\mu(dx)$  for given $\bx = \COP (x)$ fixed. For a complete mathematical formulation of conditional expectation see Appendix~\ref{sigma_algebra}. 
The conditional Helmholtz  free energy $A(\bx)$  is thus defined such that the CG probability density $\BARM(\bx)$ is of Gibbs type, i.e.
  \begin{equation}\label{PMF}
   A(\bx)=  -\frac{1}{\beta} \log \BARM(\bx) -   \frac{1}{\beta} \log Z \PERIOD
   \end{equation}

  The conditional potential of mean force (PMF)  $ \PMF(z)$ is directly related to the free energy through the reversible work theorem. \cite{Kirkwood1935, mcquarrie2000statistical}
   In many works the free energy and potential of mean force are used interchangeably.  Here we use the potential of mean force notation and write $\PMF(\bx)=A(\bx) $.   
We  define the  mean force $\FPMF: \R^{3M}\to \R^{3M} $ corresponding to the PMF defined by \VIZ{PMF}, assuming it exists,
\begin{eqnarray}\label{FPMF}
\FPMF_i (\bx) = -\nabla_{\bx_i} \PMF(\bx),\  i=1,\dots, M\PERIOD
 \end{eqnarray} 
The calculation of the potential of the mean force is a task as difficult/costly as is calculating expectations on the microscopic space. Thus one seeks   an effective potential function $\BARV(\bx)$, that approximates as well as possible the PMF,
 and  is easy to formulate and calculate.
This is the ultimate goal of all the methods  (i.e., structural-based methods, force matching, relative entropy minimization)  for systems in equilibrium, that we present in the following sections in detail. In all the above mentioned methods, one usually proposes a family of interaction potential functions  $\BARV(\bx)$ in a parametrized, $ \BARV( \bx;\theta)  $, $\theta\in  \Theta$, or a functional  form $  \BARV( \bx)$ and seeks for the optimal $\BARV^*(\bx)$ that 'best approximates' the PMF.
We denote by 
\begin{equation*}\label{CGmeasure}
\BARMV(d\bx) =\ \BARIT Z^{-1}  \exp\{ -\beta \BARV( \bx)\}d\bx\COMMA
\end{equation*}
the  equilibrium probability measure at the coarse grained configurational space for the given  CG potential function  $  \BARV( \bx)$,
where $\BARIT Z =  \int e^{ -\beta \BARV( \bx)} d\bx$ is the CG partition function.   
 
 \section{Conditional expectation and force matching}\label{CondFM}
 Force matching is based on the   observation of  a vector field $ h:\R^{3N}\to \R^{3M} $
 \begin{eqnarray}
  x \mapsto h(x) \in \R^{3M}\COMMA 
 \end{eqnarray}
from microscopic simulations or experimental observations,  and the definition of an  optimization problem in order to find an optimal estimator $G^*(\bx)$ of $h(x)$  as a function  of configurations in the coarse space $\R^{3M}$. The optimization problem 
is to find a $G^*:\R^{3M}\to\R^{3M}$ such that  the mean square error
  \begin{equation}\label{chi2}
   \Ll(G;h) =  \E_\mu\left[ \|h  - G (\COP)\|^2\right] = \int_{\R^{3N}} \|h(x) - G(\COP (x))\|^2 \mu(x)dx \COMMA
   \end{equation}
   is minimized, where  $\|\cdot\| $ denotes  the Euclidean norm in $\R^{3M}$.
In general $h(x)$  can be any observable quantity which,  eventually, for the force matching method \cite{IzVoth2005a}  that we study represents  the atomistic force on the coarse particles in configuration $x\in \R^{3N}$.
Hence, the force matching method at equilibrium is seeking for the optimal force  $G^*(\bx)$    as the  minimizer of the mean square error $\Ll(G;h)$  over a set of   proposed CG  forces $G(\bx)$ at the coarse space.

We recast the force matching optimization problem, as proposed in \cite{IzVoth2005, IzVoth2005a},  in probabilistic terms using the concept of conditional expectation and its interpretation as a projection on a subspace of observables. First, we present a well-known result in probability theory, see ref. \cite{shiryaev1996}.
We include the proof in Appendix~\ref{appendix1}  for completeness.
We  denote by $L^2(\mu) $ the space  of  mean square  integrable   vector fields with respect to  $\mu(dx)$, i.e.
 $L^2(\mu) =\{ h:\R^{3N}\to \R^{3M} | \int_{\R^{3N}} \|h(x)\|^2 \mu(dx) <\infty\}$. 
For a given CG map $\xi: R^{3N} \to R^{3M}$ we denote
$$
  L^2(\mu;\xi) = \{ g \in L^2(\mu) | \mbox{there exists $G: R^{3M}\to R^{3M}$ such that $g = G\circ \xi$}\}
$$
 the space of observables having the properties: 
  (i) Are square integrable observables with respect to the Gibbs measure $\mu(dx)$  and   (ii) are functions of the coarse variable $\COP(x)$. Property (i)  ensures the space has a geometry 
that allows an easy formulation of the concept of projections for functions, e.g. it is a Hilbert space.  
 The later property (ii) is called a  "(sub-) $\sigma$ algebra" in mathematics in the context of conditional expectations, see Appendix~\ref{AppendixCond} for further information.
\begin{lem}\label{lemmachi}
For a given $h\in L^2(\mu) $ the  minimization problem
\begin{equation*}
  \inf_{G }\Ll(G;h)= \inf_{G } \E_\mu  \left[  \left\| h - G(\xi)  \right\|^2\right] \COMMA
\end{equation*}
where $\inf$ is taken over all $G \in L^2( \mu;\COP)$ 
has the unique solution 
$$F(\bx) =\E_{\mu}[h|\xi=\bx], \ \bx \in \R^{3M}\PERIOD$$
Furthermore,
\begin{equation}\label{chirep}
  \Ll(G;h) = \Ll(F;h) +  \E_{\mu}\left[ \left\| F(\COP)-G( \COP) \right\|^2 \right] \ \text{ for any } G \in L^2( \mu;\COP) \PERIOD
\end{equation}
\end{lem}
\begin{figure}
 \includegraphics[scale=0.4]{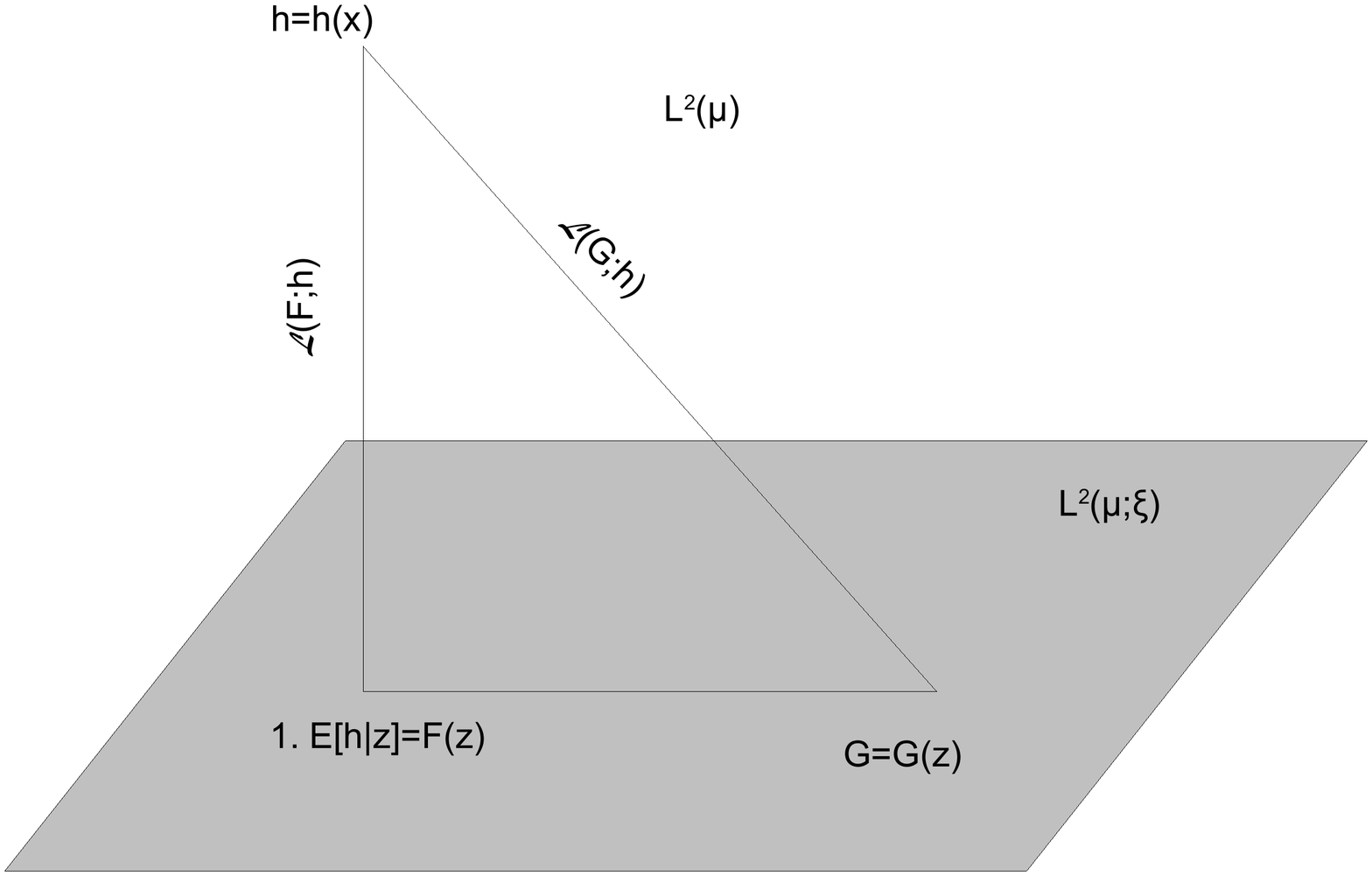}
 \caption{Projection (1)  of any microscopic   observable $h(x)$ on the space of  all CG observables $G(\COP(x))$, is given by the conditional expectation $ \E_{\mu}[h|\xi]$.   }\label{scheme1}
 \end{figure}
 A geometric description of this result is shown in Figure~\ref{scheme1}. 
In  practice the projection of $h$ is performed on a subset  of $ L^2(\mu;\COP)$, the space of feasible observables
$$\Ee  \subset  L^2(\mu;\COP)$$
that is the collection of all proposed 
CG force fields $G(\bx)$, see  Figure~\ref{scheme2}.   The set $\Ee$ may consist of   non-parametrized   or  parametrized elements, i.e., a set of splines, the span of a truncated  basis  of $ L^2(\mu;\COP)$, etc. When the minimization problem is over $ \Ee$  the solution is not necessarily  the conditional expectation $\E_{\mu}[h|\xi]$, defined by relation \VIZ{CE}, as the Lemma~\ref{lemmachi} states  since  it is possible that $ \E_{\mu}[h|\xi] \notin \Ee $, rather it is a $G^*\in \Ee$ for which relation \VIZ{chirep} holds, see the schematic in Figure~\ref{scheme2}.  In this case we say that $G^*(\bx) $, the projection of $h(x)$ on $\Ee$, is a best approximation of the  $ F(z)=\E[h|\xi]$.
 \begin{figure}
 \includegraphics[scale=0.4]{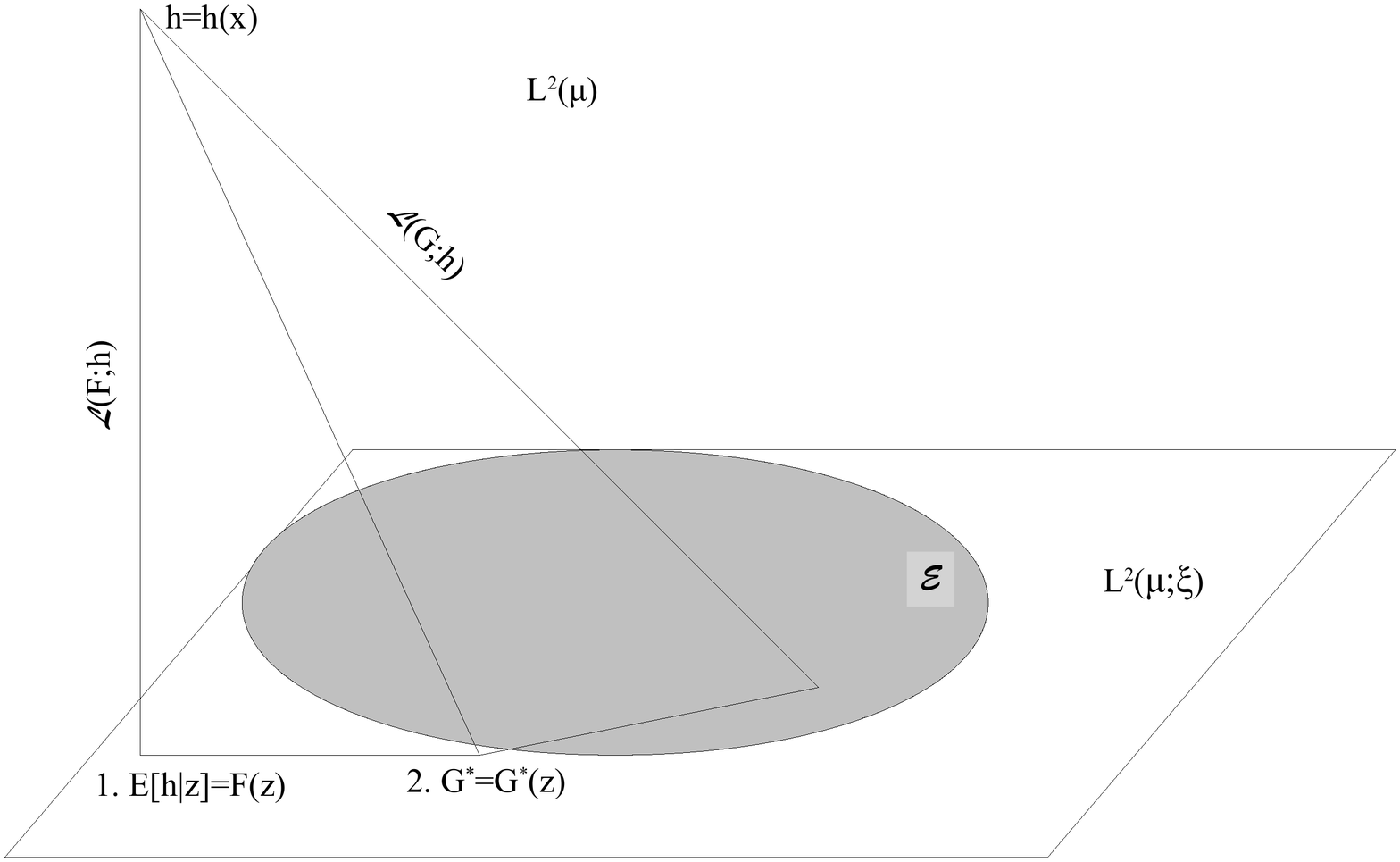}
 \caption{Geometric representation of the force matching procedure. 
 Projection (2) of the observable $h(x)$  over the set of feasible coarse observables $\Ee$. }\label{scheme2}
 \end{figure} 
 With the following theorem we state a necessary and sufficient condition that the observed quantity $h(x)$ should satisfy so that  the mean force $ \FPMF(\bx) $ \VIZ{FPMF} is best approximated with a force matching method. 
\begin{thm}[{Force matching}]\label{localmean}
Let $h\in L^2( \mu)$ such that 
\begin{equation}\label{hPMF}
\E_{\mu}[h|\bx]= \FPMF(\bx),\ \bx\in \R^{3M} \COMMA
\end{equation}
then for   $\Ee \subset L^2( \mu;\COP)$
\begin{enumerate} 
\item[a)] if  $ \FPMF\in \Ee $,
\begin{equation*}
\inf_{G\in \Ee} \Ll(G;h) = \Ll(\FPMF;h)\COMMA
\end{equation*}
\item[b)]  if $ \FPMF \notin \Ee  $   the minimizer  $G^*\in \Ee$  of $\Ll(G;h) $  satisfies
\begin{eqnarray*}
\Ll(G^*;h) =\inf_{G\in \Ee} \Ll(G;h) &=& \Ll(\FPMF;h) + \inf_{G\in \Ee}  \E_{\mu}\left[ \left\| \FPMF(\COP)-G(\COP) \right\|^2 \right] \\
&=& \Ll(\FPMF;h) +    \E_{\mu}\left[ \left\| \FPMF(\COP)-G^*(\COP )\right\|^2 \right]\PERIOD
\end{eqnarray*}
\end{enumerate}
\end{thm}
\PROOF  Is a consequence of Lemma~\ref{lemmachi}.

The result of Theorem~\ref{localmean} in combination with Lemma~\ref{lemmachi} states that \VIZ{hPMF} is a necessary and sufficient  condition for $h$ such that the corresponding force matching optimization problem has optimal, or near optimal, solution the mean force $\FPMF(\bx)$. It is thus evident that the force matching method calculates exactly the mean force if the   $\FPMF \in \Ee$. In the case that  $\FPMF \notin \Ee$ the force matching method is  best approximating the mean force in the sense that it calculates 
a $G^*$, the best approximation of $ \FPMF$ on $\Ee$, that is the element of $\Ee$  closest  to $ \FPMF$  in the $L^2$ distance, see the schematic at Figure 2. The error done   is exactly  
$$ \E_{\mu}\left[ \left\| \FPMF(\COP)-G^*(\COP) \right\|^2 \right]\PERIOD$$

The  $h(x)$ that satisfies~\VIZ{hPMF}  is called  {\it local mean force}, as in the thermodynamic integration theory. Note that  Theorem~\ref{localmean} suggests that  $h(x)$ should have a specific form   which is not at all obvious.  The purpose  of the following section is to provide closed form  representations for the  local mean force $h(x)$.

Summarizing,  the starting point, and overall goal, of a force matching method  is to  find $ \FPMF(\bx)$ for a fixed coarse graining map $\COP(x)$.  
Ideally $ \FPMF(\bx)$  could be calculated by solving a least squares problem of the form $\inf_G \E[\|\FPMF-G\|^2] $
over a   set $\Ee$ of CG models, for which the optimal solution is obviously a best approximation of $\FPMF(z)$. But although $\FPMF(z)$ is fixed it  is  not known,  thus there is a need for  using a computable quantity $h(x)$ instead of $\FPMF(\bx)$ in the minimization problem that will still has solution $\FPMF(\bx)$ or
 a best approximation  $ G^*(z)$.
 Therefore we construct an $h(x)$ that satisfies \VIZ{hPMF}, for which 
there exist many possible choices  as is proved in the following section.
The error in a force matching method has two sources, first  when the 
projection on $L^2(\mu;\COP) $ of the observed quantity $h(x)$ is not the $\FPMF(\bx)$
$$\E_{\mu}\left[  \|\FPMF(\xi) - \E_\mu[h|\xi]\|^2\right]\COMMA $$ 
and second  when the set of proposed CG forces $\Ee$ does not include $\FPMF$,
$$ \E_{\mu}\left[ \left\| \FPMF(\COP)-G^*(\COP) \right\|^2 \right]\PERIOD$$

 \section{Construction of the local mean force and systematic force matching}\label{GenFM}
 In this section we give a closed form of the local mean force $h(x)$, appearing in the force matching problem \VIZ{chi2} for which  the mean force is best approximated,
based on the statement of     Theorem~\ref{localmean}  and results from  thermodynamic integration theory (TI)   \cite{Sprik1998, denOtter, denOtterBriels, Ciccotti}.    We introduce the derived form of $h(x)$  as  the appropriate  observable   to be used in a force matching  method implementation in order to best approximate the mean force.
In  thermodynamic integration the goal is  to calculate free energy differences for a given reaction coordinate using the derivative of the free energy, see Chapter~3 in \cite{lelievre2010}. We  think of the coarse grained variable $\COP(x)$ as  a reaction coordinate, even though in the later case one does not necessarily consider coarse graining of the system.
Then we use 
the result that the derivative of free energy (the  mean force) is  given as the conditional expectation on $\COP$ of a local mean force that has a specific form, 
 a result that we state and prove    here for completeness.

 Before we state the result  we  introduce some notations and assumptions.
 We denote   $\BFD\COP(x)$ the $3M\times 3N$ matrix     with block elements $(\BFD\COP)_{ij}(x) =  \nabla_{x_j} \COP_i(x) , \ i=1,\dots,M, j=1,\dots,N$ and $ \BFJ{\COP}(x)=  \BFD\COP (x) \BFD\COP^t(x) $ the Jacobian matrix of the transformation.   For a  matrix $A$,  $A^t$   denotes its  transpose, $\text{det} A$  the  determinant and $A^{-1}$ its inverse.
 We assume that  the map $\COP$ is smooth and such   that
 $$\text{rank}\left(\BFD\COP\right) = 3M\PERIOD $$
 This assumption ensures that the Jacobian matrix of the transformation  $ \BFJ{\COP}(x)$  is non-degenerate, i.e. $\text{det} \BFJ{\COP}(x)\neq 0$ and its inverse $\BFJ{\COP}^{-1}(x)$ exists.

The assumption in Theorem~\ref{localmean}  is  that   $h(x)$ must satisfy
\begin{equation*}
\E_{\mu}[h|\bx]  = \FPMF  (\bx) \COMMA\  i=1,\dots, M\COMMA
\end{equation*}
 which, as the following Theorem states,  is not unique rather it is parametrized by a family of vector valued functions $\BFW:\R^{3N}\to\R^{3M\times 3N}$ related to the coarse graining map. 
 \begin{thm}\label{ghPMF}
 Given the CG mapping $\COP: \R^{3N} \to \R^{3M}$ and the   microscopic forces $ f_j(x)=-\nabla_{x_j} \IP(x),\ j=1,\dots,N $,  if 
 \begin{equation}\label{genhPMFW}
h_W(x) =   \BFG^{-1}_W(x) \BFW(x)  f(x) + \frac1\beta \nabla_{x} \cdot \BFG^{-1}_W(x) \BFW(x)\COMMA
 \end{equation}
 where $\BFW(x) : \R^{3M} \to \R^{3M\times 3N} $ is any  smooth   function such that
 $$ 
 \BFG_W(x)= \BFW(x) \BFD\COP^t(x)
 $$
 is invertible, then
 \begin{equation*}
 \FPMF(\bx)  = \E_{\mu}\left[h_W |  \bx\right]\PERIOD
 \end{equation*}
 \end{thm}

The above theorem  states that the choice of the local mean force $h(x)$, that is how we construct the total force for  each CG particle that corresponds to the PMF, is not  unique, nevertheless the PMF is well defined.
 For different choices of $\BFW(x)$ we can consider various force matching minimization problems,  however, the corresponding PMF is the same.    
Furthermore some of  the problems may be better than others,   simpler and cheaper to implement. At  first glance formula~\VIZ{genhPMFW} seems complicated  though  a suitable choice of $\BFW(x)$ can introduce major simplifications. 
   Note also that  in the low temperature regime, where $1/\beta \ll 1$,  term $\nabla_{x} \cdot \BFG^{-1}_W(x) \BFW(x)$ is not contributing significantly and can be  neglected.
 
 A $\BFW(x) $ always exists, at least in the case of smooth coarse graining map  that we consider, since  choosing  $\BFW(x) = \BFD\COP(x)$ we have that $ 
 \BFG_W(x)= \BFJ\COP(x) $ is invertible. In thermodynamic integration  a well t studied   choice is $\BFW(x) = \BFD\COP(x)$ \cite{Sprik1998,denOtterBriels},  that we present in the sequel as a corollary of Theorem~\ref{ghPMF}.
 \begin{cor}\label{Dcor}
  If   $\BFW(x) =\BFD\COP(x)  $ and $\textup{rank}(\BFD\COP)=3M$ then    
 \begin{equation}\label{hPMFc}
h(x) =  \BFJ{\COP}^{-1}(x) \BFD\COP(x)  f(x) +\frac1\beta \nabla_{x} \cdot  \BFJ{\COP}^{-1}(x)\BFD\COP(x)\COMMA
 \end{equation}
 where  $ \BFJ{\COP}(x) = \BFD\COP(x)\BFD\COP^t(x) $, 
 and
 $$ 
 \FPMF(\bx)  = \E_{\mu}\left[h  |  \bx\right]\PERIOD
 $$
  \end{cor}
 Note that the second term in \VIZ{hPMFc} depends on  the curvature $\nabla_{x} \cdot  \BFJ{\COP}^{-1}(x)\BFD\COP(x) $ of the sub-manifold $\Omega(\bx) = \{x: \COP(x) = \bx\}$.
The  coarse graining maps that are mainly considered in the equilibrium parametrization methods, the force matching the relative entropy minimization, are linear mappings $\COP:\R^{3N} \to \R^{3M}$
\begin{eqnarray*}
\COP_i(x) = \sum_{j=1}^{N} \zeta_{ij} x_j,\ \zeta_{ij} \in \R,\ i=1,\dots, M \COMMA
\end{eqnarray*}
for which the corresponding curvature $\nabla_{x} \cdot  \BFJ{\COP}^{-1}(x)\BFD\COP(x) =0$, since $\BFD\COP(x)=\COPL$, where $\COPL=[\zeta_{ij }\Id_3]_{i=1,\dots,M, j=1,\dots,N}$ is independent of $x$, $\Id_3$ denotes the $3\times 3$ identity matrix. The form of the local mean force is thus simplified  given in the following corollary.
 \begin{cor}\label{linearcor}
 If the CG mapping $\COP:\R^{3N} \to \R^{3M}$ is linear  with matrix $\COPL$ 
and for any  matrix    $\BFW: \R^{3M} \to \R^{3N} $  such that the matrix
 $$ 
 \BFG_W= \BFW \BFD\COP^t
 $$
 is invertible,  then for 
  \begin{equation*}\label{hWlinear}
h_W(x) = (\BFW \COPL^t)^{-1} \BFW f(x)  \COMMA
 \end{equation*} holds $$ 
 \FPMF(\bx)  = \E_{\mu}\left[h  |  \bx\right]\PERIOD
 $$
Furthermore for  $\BFW=\COPL$ and 
  \begin{equation*}
h(x) =  \BFJ{\COP}^{-1} \COPL f(x) \COMMA
 \end{equation*}
 where 
 $ \BFJ{\COP}= \COPL \COPL^t$, holds
 $$ 
 \FPMF(\bx)  = \E_{\mu}\left[h  |  \bx\right]\PERIOD
 $$

 \end{cor}
 The result of this corollary gives a compact and simpler presentation and proof of the  form for the   coarse-grained force field as described in work \cite{Voth2008a}. It provides a way to correctly calculate the total force for each CG particle from the microscopic forces corresponding to the given coarse graining map.
In a force matching method what  is  often  used  is  the  total force   acting on each CG particle as the observable quantity, i.e. for coarsening to the center of mass of K particles $h(x) = \sum_{j=1}^K  f_{j}(x)$ where $f_j(x)$  is the total force acting on particle $j$. Let us consider generally $h(x) $ of   the form $h(x)= \mathbf{B} f(x)$ for a given $3M\times 3N$ matrix $\mathbf{B}$.  The question that arises is whether with this observable we approximate the mean force associated to the specific coarse graining. In view of the result of Corollary~\ref{linearcor} the question actually is whether there exists a $\BFW$  such that 
\begin{equation*}
 \mathbf{B} f(x) = (\BFW \COPL^t)^{-1} \BFW f(x), \text{ for all } f(x)\PERIOD
 \end{equation*}
 Therefore we  are looking for a $\BFW$ such that
 \begin{equation}\label{linearW}
 \BFW ( \Id_{3N}- \COPL^t \mathbf{ B}) = O_{3M\times 3N}\COMMA
 \end{equation}
 where $O_{3M\times 3N}$ is the $3M\times 3N $ matrix with zero entries.
 The above system of equations has non-trivial solution, i.e., a non-zero matrix $\BFW$, if $\mathbf{B}$ is such that $ \Id_{3N}- \COPL^t \mathbf{B} $ is a singular matrix. 

In the following section  we study  representative  examples of molecular systems and coarse graining mappings and show in detail the application of the results of the current section, as a means of correctly calculating the CG transformation of the microscopic forces at the implementation of a force matching problem.
 
\section{Force matching formulation for linear and non-linear CG maps}
The  subject of this section is to present  analytically the  form of the local mean force $h(x)$ for  specific examples of molecular systems and  for {\it linear} and  {\it  nonlinear} CG  mappings.
Based on the result of  Theorem~\ref{ghPMF} we find  
$h(x)$   appearing  in a force matching problem,   i.e.,   when
  \begin{equation*}
  \Ll(G;h) =  \E_\mu\left[ |h  - G (\COP)|^2\right]   
   \end{equation*}
   is minimized over ${G\in\Ee}$,   and the optimal solution is a best approximation of the PMF. We consider two cases in each example,  firstly choose a $\BFW(x)$ and construct $h(x)$ and secondly we accept that the form of $h(x)$ is given and check whether satisfies the force matching condition \VIZ{genhPMFW}, i.e.  investigate whether there exists a $\BFW(x)$ appearing in \VIZ{genhPMFW}.

\subsection{ $N$-particle system under linear coarse graining maps}\label{LinearEx}
 
Let us consider a microscopic system of $N$ particles with masses $m_j,\ j=1,\dots, N$ and position vectors $x=(x_1,\dots,x_N)\in \R^{3N}$. In  the following   sections we consider different linear coarse graining maps $\xi(x)$ for which we derive explicit forms of the local mean force $h(x)$. 
Define the linear mapping $\xi: \R^{3N} \to \R^{3M}$ by
\begin{equation*}
 \xi_i(x) = \sum_{j=1}^N \zeta_{ij} x_j \in \R^{3},  \ i=1,\dots, M\COMMA
\end{equation*}
for $\zeta_{ij} \in \R$ such that $\sum_{j=1}^N{\zeta_{ij}} =1 $ for all $i=1,\dots,N$. 
The  corresponding   matrix is the $3M\times 3N $ matrix 
$\COPL=[\BFT_{ij}]_{\ontop{i=1,\dots,M}{ j=1,\dots,N}}$   where $\BFT_{ij}$ are the $3\times 3 $ blocks 
\begin{eqnarray}\label{block}
\BFT_{ij}= \zeta_{ij} \Id_{3} =
\begin{bmatrix} 
\zeta_{ij}	 & 0  				& 0  \\ 
0				 &   \zeta_{ij}	&0 \\
0 				 & 0  				& 	\zeta_{ij}
\end{bmatrix}\COMMA
\end{eqnarray}
and $\Id_3$ denotes the $3\times 3 $ identity matrix.
\subsubsection{Center of mass of N particles}
\begin{figure}
  \includegraphics[scale=0.3]{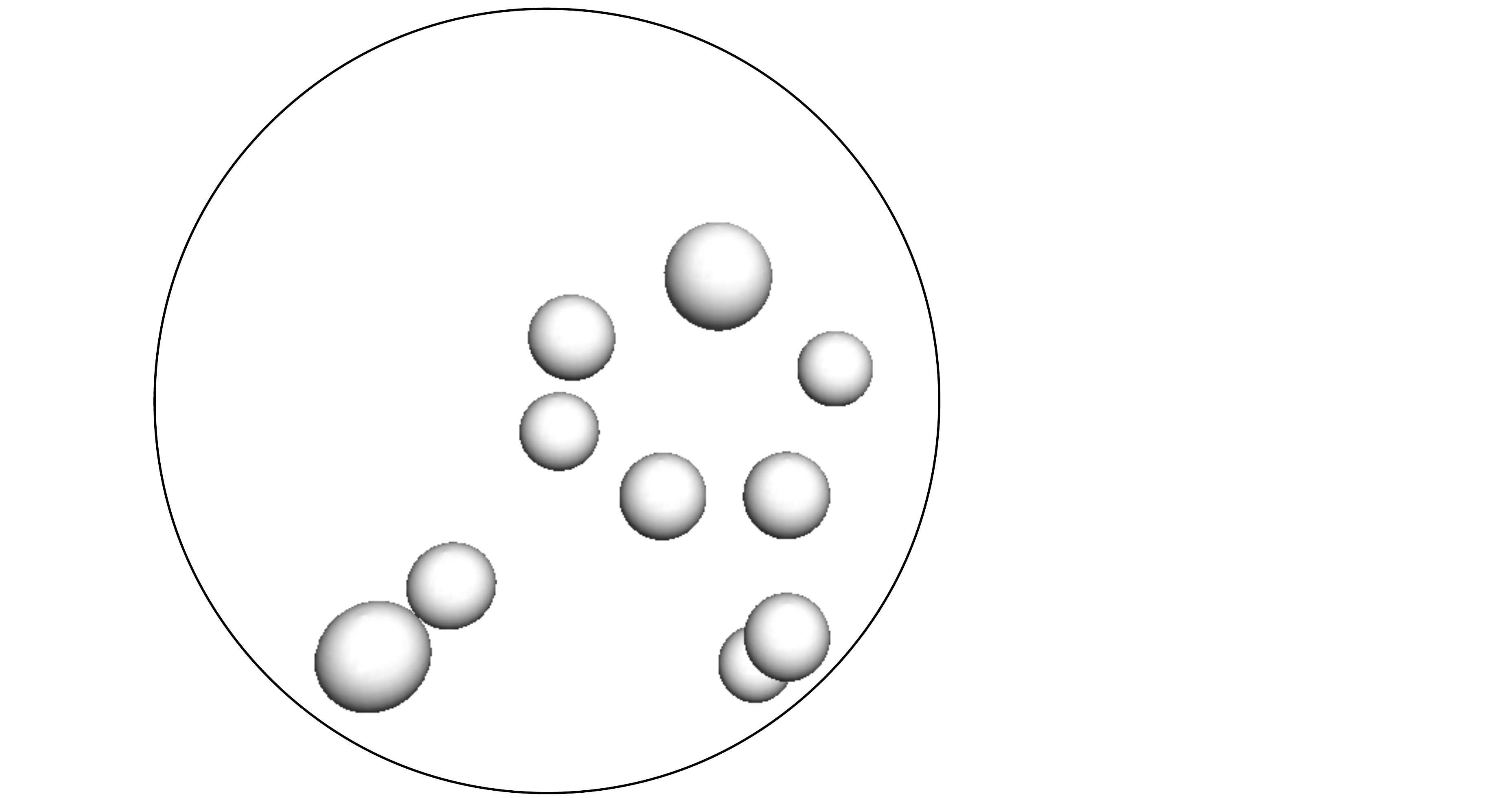}
 \caption{Coarsening a many particle system to one CG particle, the center of mass of the N particles.}\label{Center}
 \end{figure}

In this example the coarse grained variable is the center of mass of the $N$ particles, see Figure~\ref{Center}. 
That is $M=1$ and the elements of the $3\times 3N $ coarse graining mapping  matrix $\COPL$  are $\BFT_{1j}=\zeta_{1j} \Id_3,\ j=1,\dots,N$, as in \VIZ{block}, with 
$$
\zeta_{1j} = \frac{m_j}{\bm},\ j=1,\dots, N,\ \text{ where } \bm= \sum_{j=1}^N m_j  \COMMA
$$
such that  
\begin{equation*}
\xi(x) = \sum_{j=1}^N \frac{m_j}{\bm} x_j \in \R^3\PERIOD
\end{equation*}

We distinguish two cases, the first choosing a specific $\BFW$ and looking for the form of $h(x)$ and the second by choosing a local mean force $h(x)$ and looking for the existence of $\BFW$ such that the force matching indeed approximates the PMF.

\paragraph{}
  If we choose $\BFW(x) = \COPL$ then the local mean force $h(x)$ is given by 
  \begin{equation*}
h(x) =   \frac{\bm^2}{ \sum_{l=1}^N m_l^2}  \sum_{j=1}^N \zeta_{1j} f_j(x) = \frac{\bm}{ \sum_{l=1}^N m_l^2} \sum_{j=1}^N m_{j} f_j(x)\PERIOD
 \end{equation*}

 This is a result of the application of  Corollary~\ref{linearcor}, indeed we have that 
\begin{equation*}
h(x) =  \BFJ{\COP}^{-1} \COPL f(x)  \COMMA
 \end{equation*}
 where 
$$ \BFJ{\COP}= \COPL \COPL^t =\frac{1}{\bm^2} \sum_{j=1}^N m_j^2 \Id_3\PERIOD$$
Thus
\begin{equation*}
h(x) =   \frac{\bm^2}{ \sum_{l=1}^N m_l^2}  \sum_{j=1}^N \zeta_{1j} f_j(x) = \frac{\bm}{ \sum_{l=1}^N m_l^2} \sum_{j=1}^N m f_j(x)\PERIOD
 \end{equation*}
 
\paragraph{}
 We look for  $\BFW$ such that the local mean force is the {\it total  force} exerted at the center of mass  $h_{W}(x) = \sum_{j=1}^N  f_{j}(x)$. 
 Note that $h_{W}(x)$ is nonzero if external forces are present. 
 In view of relation \VIZ{linearW} such a  $3\times 3N$  $\BFW$ exists if 
   $\Id_{3N}-\COPL^t  \mathbf{B} $   is singular where $  \mathbf{B}   $ is the $3\times 3N$ matrix
$ \mathbf{B}= 
\begin{bmatrix} 
\Id_3 & \Id_3  & \dots & \Id_3  \\ 
\end{bmatrix} 
$ for which holds  $ h_{W}(x) = \mathbf{B}f(x)$.
We have that 
 \begin{eqnarray*}
 \det(\Id_{3N}-\COPL^t  \mathbf{B})= 
\begin{bmatrix} 
(1-\zeta_{11} )\Id_3&-\zeta_{11} \Id_3  & \dots & -\zeta_{11}\Id_3  \\ 
-\zeta_{12}\Id_3& (1-\zeta_{12} )\Id_3 & \dots &-\zeta_{12}\Id_3\\ 
\cdots & \dots & \dots & \cdots  \\ 
-\zeta_{1N}\Id_3& -\zeta_{1N}\Id_3  & \dots &(1-\zeta_{1N})\Id_3 \\ 
\end{bmatrix}\PERIOD
\end{eqnarray*}
 Since we assume that $\sum_j \zeta_{1j} =1$,  the sum of all column elements is zero and  $ \det(\Id_{3N}-\COPL^t  \mathbf{B})=0$. 
 Thus  there are  infinitely many nontrivial solutions   of $\BFW(\Id-\COPL^t\mathbf{B})= \mathbf{O}_{3M\times 3N} $  $\BFW=[   w_{11} \Id_3, \dots, w_{1N} \Id_3 ],$ $w_{1j}\in \R, j=1,\dots,N$. 
For example one solution is given for  $ w_{1j} =1, \  j=1,\dots,N$, that is matrix $\BFW$ in Corollary~\ref{linearcor} is 
$\BFW= 
\begin{bmatrix} 
\Id_3 & \Id_3  & \dots & \Id_3  \\ 
\end{bmatrix}\COMMA
$
for which $h_{W}(x) = \sum_{j=1}^N  f_{j}(x)$.

\subsubsection{Two CG particles coarse space}
 \begin{figure}
 \includegraphics[scale=0.4]{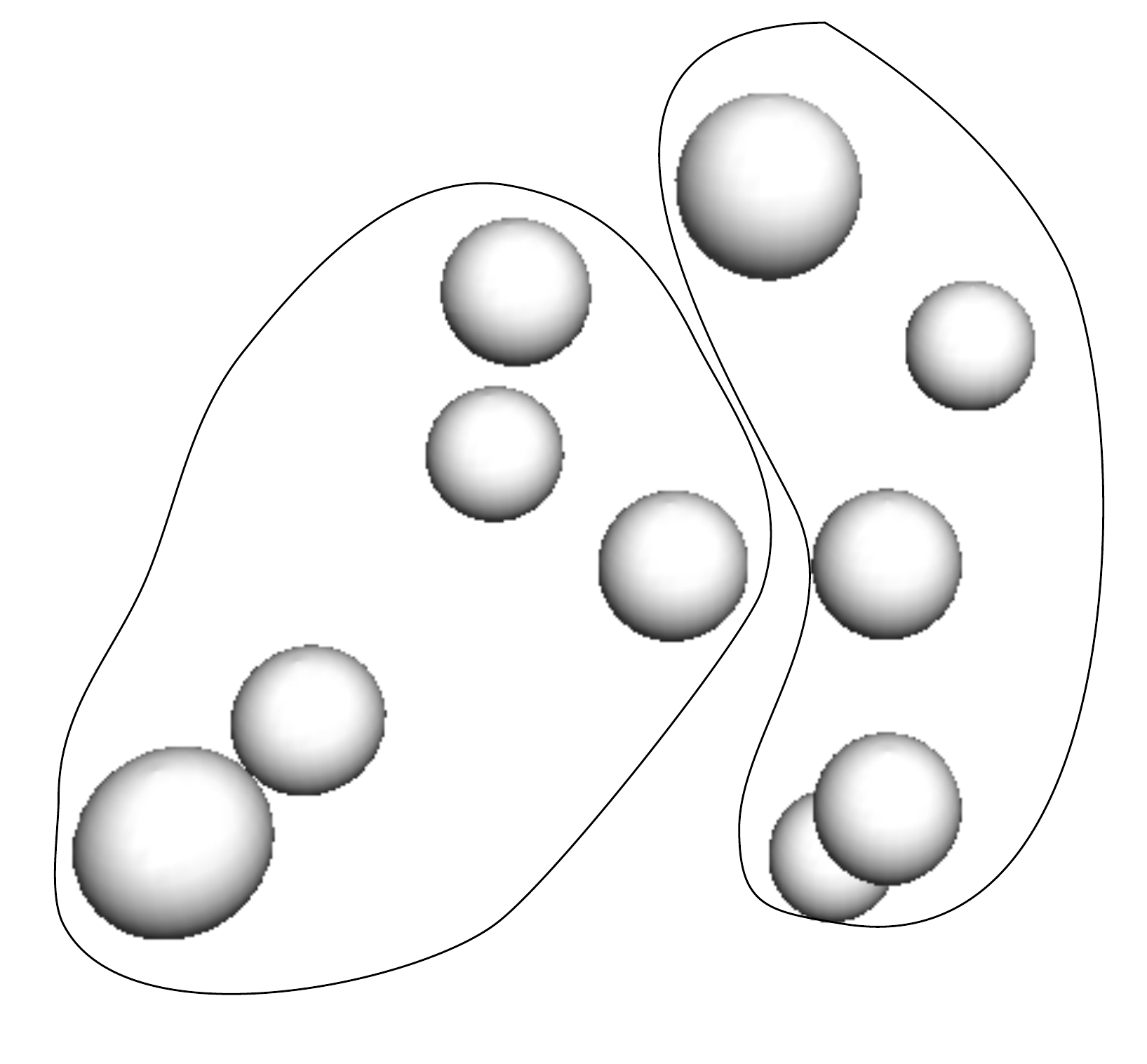}
 \caption{Coarsening a many particle system to two CG particles.}\label{TwoCG}
 \end{figure}
We consider an example where  the coarse space consists of $M=2$ CG particles, Figure~\ref{TwoCG} and  the corresponding coarse graining map is  defined by 
\begin{eqnarray*}
 \xi_1(x) = \sum_{j=1}^N \zeta_{1j} x_j,  \ \ 
  \xi_2(x) =\sum_{j=1}^N \zeta_{2j} x_j \COMMA
\end{eqnarray*}
with corresponding matrix 
\begin{eqnarray*}
\COPL= 
\begin{bmatrix} 
\BFT_{11} & \dots & \BFT_{1N} \\ 
 \BFT_{21} & \dots & \BFT_{2N}\end{bmatrix}, \  \BFT_{ij}=\zeta_{ij} \Id_3, \ i=1,2,\ j=1,\dots,N\PERIOD
\end{eqnarray*}

\paragraph{}\label{WTcase}
Let $\BFW = \COPL$, then  
\begin{equation*}
h(x) = \left(\sum_j \zeta_{1j}^2 \sum_j \zeta_{2j}^2 -\left(\sum_j \zeta_{1j}\zeta_{2j}\right)^2 \right)^{-1} 
  \begin{bmatrix} 
  \sum_j \zeta_{2j}^2  \sum_j \zeta_{1j} f_j(x) -\sum_j \zeta_{1j}\zeta_{2j} \sum_j \zeta_{2j} f_j(x)\\  
 \sum_j \zeta_{1j}^2  \sum_j \zeta_{2j} f_j(x) -\sum_j \zeta_{1j}\zeta_{2j} \sum_j \zeta_{1j} f_j(x) \end{bmatrix}  \PERIOD
  \end{equation*}
  Furthermore if  {\it each particle is contributing only to one CG particle}, 
  \begin{equation}\label{hlin2}
 h(x) =  
  \begin{bmatrix} 
   \sum_j \frac{\zeta_{1j}}{ \sum_j \zeta_{1j}^2} f_j(x)\\  
 \sum_j \frac{\zeta_{2j}}{ \sum_j \zeta_{2j}^2} f_j(x)
  \end{bmatrix}  \PERIOD
  \end{equation} 
Indeed, applying Corollary~\ref{linearcor}, we have that
\begin{eqnarray*}
 \COPL\COPL^t = 
\begin{bmatrix} \sum_j \zeta_{2j}^2 \Id_3 & \sum_j \zeta_{1j}\zeta_{2j} \Id_3 \\   \sum_j \zeta_{1j}\zeta_{2j}\Id_3  & \sum_j \zeta_{1j}^2 \Id_3\end{bmatrix}\COMMA
\end{eqnarray*}

\begin{eqnarray*}
\left(\COPL\COPL^t\right)^{-1}=\frac{1}{\text{det}(\COPL\COPL^t)}
\begin{bmatrix} \sum_j \zeta_{2j}^2 \Id_3 & -\sum_j \zeta_{1j}\zeta_{2j} \Id_3\\  -\sum_j \zeta_{1j}\zeta_{2j}\Id_3  & \sum_j \zeta_{1j}^2 \Id_3\end{bmatrix}\COMMA
\end{eqnarray*}
where
$ \text{det}(\COPL\COPL^t) = \sum_j \zeta_{1j}^2 \sum_j \zeta_{2j}^2 - (\sum_j \zeta_{1j}\zeta_{2j})^2$, 
and   the local mean force is given by 
\begin{equation*}
h(x) = \frac{1}{\text{det}(\COPL\COPL^t)} 
  \begin{bmatrix} 
  \sum_j \zeta_{2j}^2  \sum_j \zeta_{1j} f_j(x) -\sum_j \zeta_{1j}\zeta_{2j} \sum_j \zeta_{2j} f_j(x)\\  
 \sum_j \zeta_{1j}^2  \sum_j \zeta_{2j} f_j(x) -\sum_j \zeta_{1j}\zeta_{2j} \sum_j \zeta_{1j} f_j(x) \end{bmatrix}  \PERIOD
  \end{equation*} 
 If we consider that {\it each particle is contributing only to one CG particle},
i.e. $\zeta_{1j} \zeta_{2j} = 0 $, the form of $h(x) $ is  simplified,  since 
\begin{equation*}
\left(\COPL\COPL^t\right)^{-1} =  
  \begin{bmatrix} 
    \frac{ 1 }{ \sum_j \zeta_{1j}^2} & 0\\  
 0 &     \frac{ 1}{ \sum_j \zeta_{2j}^2}\end{bmatrix}  \COMMA
  \end{equation*}
  and becomes
 \begin{equation*} 
 h(x) =  
  \begin{bmatrix} 
   \sum_j \frac{\zeta_{1j}}{ \sum_j \zeta_{1j}^2} f_j(x)\\  
 \sum_j \frac{\zeta_{2j}}{ \sum_j \zeta_{2j}^2} f_j(x)
  \end{bmatrix}  \PERIOD
  \end{equation*} 
Note that when $\zeta_{ij} = m_j/\sum\limits_{k\in C_i  } m_k, \ i=1,2$, and $C_i=\{j: \zeta_{ij}\neq 0\}$, or equivalently $C_i=\{  j: \text{ particle } j  \text{ contributes to CG particle } i\}$, then from relation~\VIZ{hlin2} we have
  \begin{equation*}
h(x) =  
  \begin{bmatrix} 
   \sum_{j\in C_1} \frac{\bm_1 m_j}{   \sum_{k\in C_1} m_k^2} f_j(x)\\  
 \sum_{j\in C_2} \frac{\bm_2 m_j}{   \sum_{k\in C_2} m_k^2} f_j(x)
  \end{bmatrix}  \COMMA
  \end{equation*} 
  denoting $\bm_i =  \sum_{j\in C_i}  m_j, \ i=1,2$,
thus
   \begin{equation*}
h(x) =  
  \begin{bmatrix} 
   \sum_{j\in C_1}   f_j(x)\\  
 \sum_{j\in C_2}   f_j(x)
  \end{bmatrix}  \COMMA
  \end{equation*}    
 when all particles have equal mass $m_j=m,\ j=1,\dots,N$.
 
 \paragraph{}\label{total}
In this case we look for $\BFW$ such that  the local mean force is  
 \begin{equation*}
h_W(x) =    \begin{bmatrix} 
   \sum_{j\in C_1} f_j(x)\\  
  \sum_{j\in C_2}  f_j(x)
  \end{bmatrix} \COMMA C_i=\{j: \zeta_{ij}\neq 0\}, i=1,2 \PERIOD
  \end{equation*} 
We show that such a $\BFW$ exists if each particle is contributing only to one CG particle, that is $ \zeta_{1j}\zeta_{2j}=0$ for all $j=1,\dots,N$. 
Indeed, following relation \VIZ{linearW} we write 
\begin{equation*}
h_W(x) =    \mathbf{B}  f(x),\quad \mathbf{B} =\begin{bmatrix} 
   \delta_{\zeta_{11}} \Id_3 \dots  \delta_{\zeta_{1N}} \Id_3 \\  
   \delta_{\zeta_{21}} \Id_3 \dots  \delta_{\zeta_{2N}} \Id_3
  \end{bmatrix}  \COMMA
  \end{equation*} 
where  $\delta_{\zeta_{ij}} = 1  $ if $\zeta_{ij} \neq 0$ and  $ \delta_{\zeta_{ij}} = 0  $ if $\zeta_{ij} = 0$ and 
 calculate the $\det(\Id_{3N} - T^t B) $, where 
 \begin{equation*}
\!\Id_{3N}-\COPL^t \mathbf{B} \!=\!\!\begin{bmatrix} 
  (1-\zeta_{11} \delta_{\zeta_{11}} - \zeta_{21} \delta_{\zeta_{21}}) \Id_3& -(\zeta_{11} \delta_{\zeta_{12}} + \zeta_{21} \delta_{\zeta_{22}})\Id_3 & \dots & -(\zeta_{11} \delta_{\zeta_{1N}} + \zeta_{21} \delta_{\zeta_{2N}})\Id_3\\  
   \vdots & \vdots &\vdots&\vdots\\
  -(\zeta_{1N} \delta_{\zeta_{11}} + \zeta_{2N} \delta_{\zeta_{21}})\Id_3 & -(\zeta_{1N} \delta_{\zeta_{12}} + \zeta_{2N} \delta_{\zeta_{22}})\Id_3   & \dots  &(1-\zeta_{1N} \delta_{\zeta_{1N}} - \zeta_{2N} \delta_{\zeta_{2N}})\Id_3\\
  \end{bmatrix} \!\!\! \PERIOD
  \end{equation*} 
Using the assumption that $\sum_{j\in C_i} \zeta_{ij} =1, i=1,2$ and properties of matrix determinants we have that
\small{
  \begin{equation*}
\!\!\det(\Id_{3N}-\COPL^t \mathbf{B} )\!\!=\!\!\det\begin{bmatrix} 
  (1-  \delta_{\zeta_{11}} -   \delta_{\zeta_{21}}) \Id_3& (1- \delta_{\zeta_{12}} -   \delta_{\zeta_{22}})\Id_3 & \dots & (1- \delta_{\zeta_{1N}} -   \delta_{\zeta_{2N}})\Id_3\\  
   \vdots & \vdots &\vdots&\vdots\\
 \! -(\zeta_{1N} \delta_{\zeta_{11}}\! + \!\zeta_{2N} \delta_{\zeta_{21}})\Id_3 & \ -(\zeta_{1N} \delta_{\zeta_{12}} \!+\! \zeta_{2N} \delta_{\zeta_{22}})\Id_3   & \dots  &(1-\!\zeta_{1N} \delta_{\zeta_{1N}}\! -\! \zeta_{2N} \delta_{\zeta_{2N}})\Id_3\!\\
  \end{bmatrix} \!\!\! \PERIOD
  \end{equation*} 
  }
 Thus we see that $ \det(\Id_{3N}-\COPL^t \mathbf{B} ) = 0 $ if $1-  \delta_{\zeta_{1j}} -   \delta_{\zeta_{2j}} =0 $ for all $j=1,\dots,N $, i.e., when $ \zeta_{1j}\zeta_{2j}=0 $ that means particle $j$ {\it contributes only to one CG particle}.
 
 Assume  now that there exist only one particle $k$ that contributes to both CG particles. Let for simplicity choose $k=1$, and $\zeta_{21}=1, \zeta_{2j}=0, \ j\neq 1 $, and  $\zeta_{1j}\neq 0  $ for all $ j=1,\dots, N$ s.t. $\sum_j\zeta_{1j}=1 $ then  
 \begin{equation*}
\det(\Id_{3N}-\COPL^t \mathbf{B} )=\begin{bmatrix} 
  -  \Id_3& 0_3 & \dots & 0_3\\  
    - \zeta_{12}  \Id_3 & (1- \zeta_{12} ) \Id_3   & \dots  & -\zeta_{12}  \Id_3\\
   \cdots & \dots &\dots &\dots\\
  - \zeta_{1N}  \Id_3 & - \zeta_{1N}  \Id_3   & \dots  &(1-\zeta_{1N}  )\Id_3\\
  \end{bmatrix}  \COMMA
  \end{equation*} 
  which is nonzero.  Thus we found an example of CG map for which there does not exist any  $\BFW$ s.t. $h_W(x) =    \begin{bmatrix} 
   \sum_{j\in C_1} f_j(x)\\  
  \sum_{j\in C_2}  f_j(x)
  \end{bmatrix} \COMMA C_i=\{j: \zeta_{ij}\neq 0\} \PERIOD $ 
  This suggests that for this CG map one should choose a $\BFW$ and then construct the $h(x)$ in order to achieve the PMF approximation with the force matching, as is calculated for example in case~\ref{WTcase}.
  
  \subsubsection{Two particles contributing to each CG particle}
 \begin{figure}
  \includegraphics[scale=0.4]{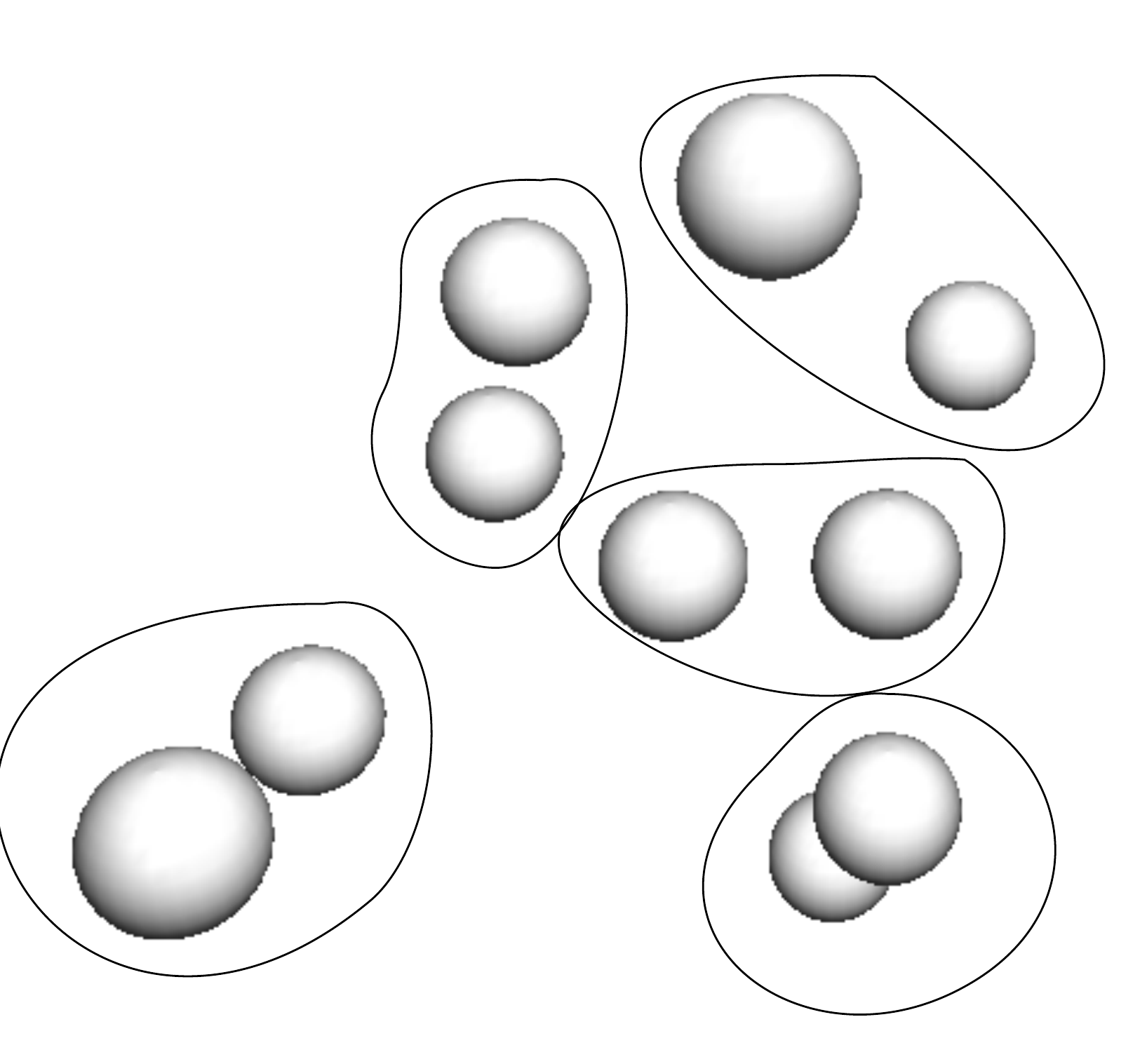}
 \caption{Coarsening a many particle system with two particles per CG particle.}\label{Pairs}
 \end{figure}
With this example we examine the coarse graining where each CG particle   is the average  of two    particles position vectors  which  contribute only to that   CG particle, Figure~\ref{Pairs}. That is, assuming that the number of particles $N$ is even, the number of CG particles is $M =N/2$  and  the mapping is defined by
\begin{eqnarray*}
 &&\xi_i(x) =  \zeta_{i,2i-1} x_{2i-1} + \zeta_{i,2i} x_{2i}, \ i=1,\dots,M\COMMA
\end{eqnarray*}
for $\zeta_{ij}\in \R$, and $\zeta_{ij}= 0$ if $ j\neq 2i-1, 2i$, such that $ \zeta_{i,2i-1}   + \zeta_{i,2i}  =1 $ for all $i=1,\dots,M$.
The $3M\times 3N$ matrix of the linear mapping $\xi$ is 
\begin{eqnarray*}
\COPL= 
\begin{bmatrix} 
\BFT_{11} & \BFT_{12} & 0 & 0 &0 & \dots & \dots & 0 &0 \\ 
0				 &   0 			  &\BFT_{23} & \BFT_{24} & 0 & \dots & \dots & 0 &0\\
 \vdots                &  			&	 & 	&	&\BFT_{i,2i-1} & \BFT_{i,2i}  	&0	& \vdots\\
0 				 & 0  		&0&	   \dots    &\dots &\dots   &0 & \BFT_{MN-1} & \BFT_{MN}  
\end{bmatrix}\COMMA
\end{eqnarray*}
where $ \BFT_{ij} = \zeta_{ij} \Id_3,\ i=1,\dots, M, \ j=1, \dots, 2M$. 

\paragraph{}
 Let $\BFW=\COPL$, applying Corollary~\ref{linearcor} we have that 
\begin{equation*}
h_i(x) = \frac{1}{\zeta_{i,2i-1}^2   + \zeta_{i,2i}^2 } \left( \zeta_{i,2i-1} f_{2i-1}(x) + \zeta_{i,2i} f_{2i}(x)\right), \ i=1,\dots,M \PERIOD
  \end{equation*} 
Let $m_j$ denote the mass of the $j-{th}$ particle and set 
 $$\zeta_{i,2i-1} = \frac{m_{2i-1}}{m_{2i-1}+m_{2i }}, \  \zeta_{i,2i}=\frac{m_{2i }}{m_{2i-1}+m_{2i }}   \   \ i=1,\dots,M\COMMA$$
  then if   $ m_{2i-1} = m_{2i }, \ i=1,\dots,M$ we can have 
  \begin{equation*}
h_i(x)  =  f_{2i-1}(x) +  f_{2i}(x), \ i=1,\dots,M \PERIOD
  \end{equation*} 
 
 \paragraph{}
 We show   that there exists a family of $3M\times 3N$ matrices $\BFW$ appearing in Theorem~\ref{ghPMF}, such that 
  \begin{equation*}
h_{W,i}(x) =  f_{2i-1}(x) +  f_{2i}(x), \ i=1,\dots,M \PERIOD
  \end{equation*} 

 Let $\BFW$ with block entries $ w_{ij} \Id_3, i=1,2,\dots, M \  j=1,\dots,N $.  The later equality holds if,  in view of \VIZ{linearW},  
$ \BFW ( \Id_{3N}- \COPL^t \mathbf{ B}) = O_{3M\times 3N}$
  where 
  \begin{eqnarray*}
\mathbf{B} = 
\begin{bmatrix} 
\Id_3 & \Id_3 & 0 & 0 &0 & \dots & \dots & 0 &0 \\ 
0				 &   0 			  &\Id_3 & \Id_3& 0 & \dots & \dots & 0 &0\\
 \vdots                &  			&	 & 	&	&\Id_3 & \Id_3  	&0	& \vdots\\
0 				 & 0  		&0&	   \dots    &\dots &\dots   &0 &\Id_3& \Id_3 
\end{bmatrix}\PERIOD
\end{eqnarray*}
We have that $\det(\Id_{3N}- \COPL^t \mathbf{B}  ) =0$ since $ \zeta_{i,2i-1}   + \zeta_{i,2i}  =1 $ for all $i=1,\dots,M$. 
Therefore, according to \VIZ{linearW}, there exist infinitely many  $M\times (2M)$ matrices $   \BFW $   such that $ \BFW(\Id_{3N}-\COPL^t \mathbf{B} ) =0$, that gives 
 \begin{equation*}
  w_{i,2j}=w_{i,2j-1},\ \text{ for all } i,j =1,\dots,M\COMMA
  \end{equation*}                     
and
  \begin{equation*}
\left(h_W(x) \right)_i=  f_{2i-1}(x) +  f_{2i}(x), \ i=1,\dots,M \PERIOD
  \end{equation*} 
                     
Note that in all the above examples with linear CG maps the form of $h(x)$ can also be written in the form   
\begin{equation*}
h(x) =    \sum_{j=1}^N\frac{d_{ij}}{ \zeta_{ij}} f_j(x), i=1,\dots,M \COMMA
 \end{equation*}
when appropriately choosing constants  $ d_{ij}$. This fact proves that our approach reproduces the results in  Noid et.al \cite{Voth2008a} and  
is indeed an extension that holds for any nonlinear CG map, that we show with examples in the following section.                                            
\subsection{Force matching formulation and non-linear CG maps}\label{ReactionC}
In this section we examine the application of the force matching method with   examples where we consider that the coarse graining mapping corresponds to  a reaction coordinate, that is in principle   a   nonlinear mapping $\COP:\R^{3N}\to \R^m$. 

 We borrow the   example  from ref.
\cite{denOtter}, where   the corresponding free energy differences and PMF were calculated explicitly using generalized coordinates. 
Here we  only consider the mapping to the reaction coordinate and a proper selection of $\BFW(x)$ appearing in \VIZ{genhPMFW}, as is also remarked in  \cite{chipot2007free} Section 4.4.
In this example the microscopic model    is  a single   molecule consisting of three atoms. Let $x_j\in \R^3$, $j=1,2,3$ denote the position vectors of the atoms, see Figure~6. 

\subsubsection{Bending angle}
The coarse variable  is the {\it bending} angle
\begin{figure}\label{threeatom}
 \includegraphics[scale=0.4]{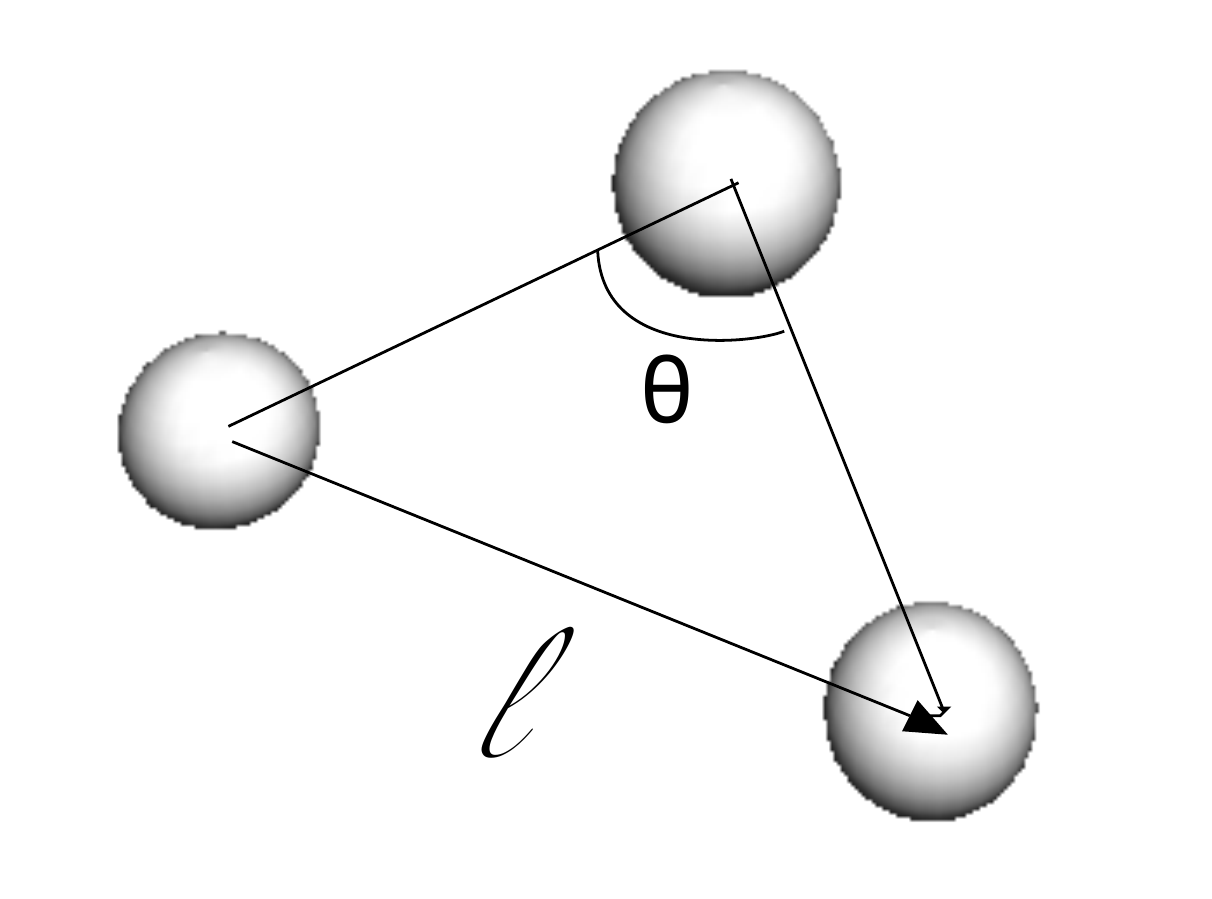}
 \caption{Three atom molecule a) bending angle $\theta$ b) end-to-end vector $\vec{\ell}$.}
 \end{figure}
 $\theta=<x_1x_2x_3$   , see Figure~6,
\begin{equation}\label{bending}
 \xi:\R^{9} \to (0, 2\pi),\quad
  \xi(x) = \text{acos}\frac{<x_3-x_2,x_1-x_2>}{\|x_3-x_2\|\|x_1-x_2\|}:=\theta\COMMA
\end{equation}
where $<\cdot,\cdot>$, $\|\cdot\|$ denote the Euclidean inner product and norm in $\R^{3}$ respectively.
Applying Corollary~\ref{Dcor}, that is choosing $\BFW(x) = \BFD\COP(x) $,  the local mean force is 
 \begin{equation*} 
h(x) =  \BFJ{\COP}^{-1}(x) \BFD\COP(x)  f(x) +\frac1\beta \nabla_{x} \cdot  \BFJ{\COP}^{-1}(x)\BFD\COP(x)\COMMA
 \end{equation*}
 where  $ \BFJ{\COP}(x) = \BFD\COP(x)\BFD\COP^t(x) $.
Here  $ \BFD\COP(x) = \left(\nabla_{x_1} \xi(x), \nabla_{x_2}\xi(x), \nabla_{x_3} \xi(x)\right) \in \R^{1\times 9}$ 
where
$$\nabla_{x_j} \xi(x) = -\frac{1}{\sin(\xi(x))} \nabla_{x_j} \frac{<x_3-x_2,x_1-x_2>}{\|x_3-x_2\|\|x_1-x_2\|}, \ j=1,2,3\COMMA$$
and 
   $ \BFJ{\COP}(x) = \|  \BFD\COP(x) \|^2 \in \R$.
   Thus 
    \begin{equation*} 
h(x) =  \frac{1}{\|  \BFD\COP(x) \|^2}\BFD\COP(x)  f(x) +\frac1\beta \nabla_{x} \cdot \left( \frac{1}{\|  \BFD\COP(x) \|^2}  \BFD\COP(x)\right)\COMMA
 \end{equation*}
 where $\COP(x)$ is given by \VIZ{bending}.

  \subsubsection{End to end distance}
 Let us now choose the end to end  distance  $\|\ell_{13}\|$  as a coarse variable,
 \begin{equation}\label{eted}
 \xi :\R^{9} \to (0,\infty), \quad \xi(x) = \|x_1-x_3\| :=\|\ell_{13}\|\COMMA
\end{equation}
for which 
$ \BFD\COP(x) $=
$ \|x_1-x_3\|^{-1}\left(x_1-x_3,   0_3, x_3-x_1\right) $ and $ \BFJ\COP= \BFD\COP \BFD\COP^t =2$.
 Applying Corollary~\ref{Dcor}, we have
     \begin{eqnarray*}
h(x) &=& \frac{1}{2}  \BFD\COP(x)  f(x) +\frac1\beta \nabla_{x} \cdot \left( \frac{1}{2}  \BFD\COP(x)\right)\\
&=&  \frac{<x_3-x_1,  f_3(x)- f_1(x)>}{2 \|x_3-x_1\|}    +\frac3\beta  \COMMA\\
     \end{eqnarray*} 
since 
 $\nabla_{x} \cdot   \BFD\COP(x) = 6 $, following the definition \VIZ{eted}.
 
 \noindent
 {\bf Remark.} A coarse variable that is of interest in molecular systems is the end-to-end vector     $  \ell_{13} = x_3-x_1\in \R^3$,   the corresponding map is linear  with  
  \begin{equation*} 
 \xi :\R^{9} \to \R^{3}, \quad \xi(x) =  x_3-x_1 :=\ell_{13}\COMMA
\end{equation*} 
 The  mapping has the $ 3\times 9  $ corresponding matrix
\begin{equation*}
\COPL= \BFD\COP=
\begin{bmatrix} 
-1 & 0 &0 &0 &0 &0&1&0&0 \\ 
 0 & -1 &0 &0 &0 &0&0&1&0\\
 0 & 0 &-1 &0 &0 &0&0&0&1
 \end{bmatrix}\PERIOD
 \end{equation*}
 Since the mapping is linear we can apply Corollary~\ref{linearcor}, which  gives  
     \begin{equation*} 
h(x) =  \frac{1}{2}\BFD\COP(x)  f(x) = \frac12 \left( f_3(x) - f_1(x)\right)\PERIOD
 \end{equation*}
\section{Force matching and information-based projections }\label{REandFM}
In this section we show there is a strong link between coarse-graining viewed as minimization of relative entropy and CG derived from force matching optimization principle in $L^2(\mu)$ presented in Sections~\ref{CondFM} and \ref{GenFM}.
We first start the discussion with a brief outline of the relative entropy minimization and continue with its the relation with force matching. Finally, we include  a brief description  of structural based methods   in order to provide a complete view of the  methods for potential of mean force approximations in coarse graining. 
  \subsection{Relative entropy}
 The relative entropy approach \cite{Shell2008, Shell2011}  considers the minimization of the relative entropy functional 
\begin{equation}\label{REfunc}
 \min_{\BARV\in \Vv}\RELENT{\mu}{ \mu_{\BARV}}=   \min_{\BARV\in \Vv}\E_{\mu}\left[ \log\frac{\mu}{\mu_{\BARV}}\right]
\end{equation}
over a  space $  \Vv =\{\BARV | \BARV :\R^{3M} \to \R\}$ of interaction potentials. If the CG potentials are parametrized with $\theta \in \Theta$  the minimization is considered over  the parameter space $  \Theta$. The minimization problem  is based on the properties of the relative entropy  
a) $ \RELENT{\mu}{\pi}\ge 0$ for all probability measures $\mu,\pi$ and b)$ \RELENT{\mu}{\pi} = 0$ if and only $\mu \equiv \pi$.

The relative entropy  $ \RELENT{\mu}{ \mu_{\BARV}}$ is a pseudo-distance between the microscopic Gibbs measure $\mu(x)\propto e^{-\beta \IP(x)} dx$ and a back-mapping of the proposed Gibbs measure at the CG space $\BARMV(\bx) \propto e^{-\beta \BARV(\bx)}d\bx$
$$ 
\mu_{\BARV}(x)=\BARMV(\bx)\nu(x|\bx)\COMMA
$$
 associated with the proposed interaction potential ${\BARV}(\bx)$ where
\begin{equation}\label{recon}
\nu(x|\bx)=\frac{\gamma(x)}{Z_{\gamma}( \bx)}\ \text{ with } \gamma\text{    any nonnegative } L^1(\R^{3N}) \text{  function, }
\end{equation}
and 
$$ \  Z_{\gamma}(\bx) = \EXPECT[{\gamma}|{\bx}]=\int_{\Omega(\bx)} \gamma(x) dx, \   \int_{\Omega(\bx)} \nu(x|\bx) dx =1 \PERIOD $$ 
Recall that $\Omega(\bx) =  \{ x\in \R^{3N}: \COP (x)= \bx\}$. 
The measure $\nu(x|\bx)$ is  a normalized conditional probability of sampling an atomistic configuration $x$ given a CG configuration $\bx$ (microscopic reconstruction).  A  mathematical  formulation of microscopic reconstruction is presented in our  work \cite{KPR} while probabilistic reconstruction methodologies are proposed and tested in
\cite{KPR,TT,KKPV, MulPlat2002,tsop,Harmandaris2006a}.

 The   difference  in relative entropy between $\mu(x)$ and $\mu_{\BARV}(x)$ is written  
\begin{equation}\label{REerror}
\RELENT{\mu}{ \mu_{ \BARV}} = \RELENT{\BARM}{ \BARM_{ \BARV}} + \int \RELENT{\mu(\cdot|z)}{ \nu(\cdot|z) }\BARM(dz)
\end{equation}
where  $\BARM(z)=\E[\mu|z]$    is the exact coarse grained measure \VIZ{exactCGm}, and  $\mu(x|z)$ is the unique measure $\frac{d\mu}{d\BARM}(x|z) $, i.e., such that $\mu(dx) = \BARM(dz) \mu(dx|z) $. 
Relation \VIZ{REerror} shows that the  difference  is composed from two parts a) the error in the approximation of the exact Gibbs measure $ \BARM(z)$  corresponding to the $\PMF(z) $ by $\BARM_{ \BARV}$, $\RELENT{\BARM}{ \BARM_{ \BARV}}$, and  b) the error in  reconstruction, $\int \RELENT{\mu(\cdot|z)}{ \nu(\cdot|z) }\BARM(dz)$, that is the error in approximating $ \mu(x|z)$ by $\nu(x|z)$. 

In the relative entropy minimization method, as defined by Shell et.al. \cite{Shell2008, Shell2009}  $\gamma(x)=1$  assigning the same probability to all atomistic configurations $x$ that map to the same $\bx$.
The reconstruction measure is  the uniform distribution 
 $ \nu(x|z)= 1/|\Omega(z)|$, where $|\Omega(z)|$ is the volume of the set $\Omega(z)$, and the error introduced is 
 $$ 
 \RELENT{\mu(\cdot|z)}{ \nu(\cdot|z) } = \log |\Omega| + \int \mu(x|z) \log\mu(x|z) dx|z\PERIOD
 $$
 Note that for this choice of reconstruction  the error does not depend on the proposed approximating potential $\BARV(\bx)$, the error  is constant for any   $\BARIT \IP(\bx)\in \Vv$. 
In the ideal case where $ \gamma(x)=\mu(x)$ the reconstruction is considered exact,  there is no reconstruction error since $ \nu(x|z)  = \mu(x|z)$ and
$\RELENT{\mu(\cdot|z)}{ \nu(\cdot|z) } = 0$, and the minimization problem is equivalent to 
$\min_{\BARV\in \Vv}\RELENT{\BARM}{ \BARM_{ \BARV}}.$

In view of the last two observations  it is verified  that the relative entropy minimization method, with uniform or exact reconstruction, is indeed approximating the potential of mean force $\PMF(\bx)$ since the minimization problem   $\min_{\BARV\in \Vv}\RELENT{\mu}{ \mu_{ \BARV}}$ is equivalent to the 
$$\min_{\BARV\in \Vv}\RELENT{\BARM}{ \BARM_{ \BARV}}\PERIOD$$
 
\subsection{Relative entropy and Force matching}
The goal  of the last part of this section is to compare the  force matching method with the relative entropy minimization method.
The common point of both methods is their relation to the PMF.  The relative entropy is directly related with the PMF through relation \VIZ{REerror} while  the force matching method at  equilibrium   approximates the PMF  if, as stated in Theorem~\ref{ghPMF}, the local mean force  $h(x)$   is such that  $\FPMF(\bx)  = \E_{\mu}\left[h |  \bx\right]  $.

As discussed in the previous section,  a reasonable choice for the reconstruction   is   $\gamma(x)=\mu(x)$, the equilibrium Gibbs measure \cite{Noid2011}, thus 
$$Z_{\gamma}(\bx) =  \BARM(\bx) \PERIOD$$ 
Practically this choice of $\gamma(x)$  means that we sample from the Gibbs measure using constraints on $\bx$. One can easily check  that the relative entropy  $ \RELENT{\mu}{ \mu_{ \BARV}}$  
for $\gamma(x) =\mu(x) $   is rewritten as
\begin{equation*} 
 \RELENT{\mu}{ \mu_{ \BARV} } = 
  \E_{\BARM}\left[ \log \frac{\BARM(\bx)}{\BARMV (\bx)}\right]=\RELENT{\BARM}{ \BARM_{ \BARV}}\PERIOD
  \end{equation*}
Based on the above equality   and  the properties of the relative entropy  
we can see that the minimum value of  $ \RELENT{\mu}{\mu_{ \BARV} } $  is given when 
 $ \BARM_{ \BARV}(z) = \BARM(\bx)$ corresponding to the PMF $\PMF(\bx)$,  under the assumption that the 
reconstruction probability $\nu(x|\bx)$ is exact, i.e., $\RELENT{\mu(\cdot|z)}{ \nu(\cdot|z) } = 0$. 
 
 With the following theorem we compare the relative entropy minimization and the force matching methods   under the assumptions that both approximate the PMF $\PMF$, in the sense discussed in the previous sections, i.e.  $\FPMF(\bx)  = \E_{\mu}\left[h |  \bx\right]  $ in force matching and  $\gamma(x)=\mu(x)$ in relative entropy minimization.  
 
\begin{thm}[{Relative entropy and force matching at equilibrium}]\label{comREFM}
\mbox{}

\noindent
Consider a microscopic system in $ \R^{3N}$ at equilibrium, characterized by the interaction potential $\IP(x)$ and the Gibbs measure $\mu(x)$.
 Let $\COP:\R^{3N} \to \R^{3M} $ be a CG mapping, $ \BARIT \IP(\bx) \in \Vv$ be a family of  interaction potentials  on the coarse space $ \R^{3M}$ with  Gibbs measure $\BARMV(d\bx)$  and $h:\R^{3N}\to \R^{3M}$   such that
 $\FPMF(\bx)  = \E_{\mu}\left[h |  \bx\right]  $. Let  $\mu_{ \BARV}(x)=\BARMV(\bx)\mu(x|\bx)$ where $ \mu(x|\bx) = \frac{d\mu}{d\BARM}(x|z)$.
Consider the following  two minimization problems at equilibrium
$$\min_{ \BARV \in \Vv}  \RELENT{\mu}{  \mu_{ \BARV} }  \COMMA
\ \text{(Relative entropy),}$$
and 
$$\min_{G\in \Ee} \Ll(G;h) = \min_{G\in \Ee} \E_\mu  \left[  \left\| h - G (\COP)  \right\|^2\right], \ \text{(Force matching),} $$
where  $ \Vv=\{ \BARV:\R^{3M} \to \R\}$, $\Ee=\{ G\in L^2 \text{ s.t. }  G_i(\bx) =  -\nabla \BARIT \IP(\bx),  \BARV\in \Vv\}$. 

 Then the leading term at the relative entropy approach  is the square of the potential difference
 \begin{equation}\label{RERphi} 
  \min\limits_{ \BARV \in \Vv}  \RELENT{\mu}{\mu_{ \BARV} } = \beta^2 \min\limits_{ \BARV \in \Vv} 
  \E_{\BARM}\left[    \left( \PMF  -   \BARV \right) ^2 \right]  + \BIGO\left(\beta^3
  \E_{\BARM}\left[ | \BARV - \PMF|^3\right]  \right)  \COMMA
  \end{equation} 
  where $\BIGO(g) $ denotes a quantity  bounded by $g$,
and the   force matching minimizes the square  of  the potential gradients difference
\begin{equation}\label{chiphi}   
 \min_{G\in \Ee} \Ll(G;h) =  \Ll(\FPMF;h) + \beta^2 \min_{ \BARV \in \Vv}  \E_{\BARM }\left[ \|  \nabla  \left(  \PMF  -   \BARV  \right) \|^2  \right]   \PERIOD
\end{equation}
  \end{thm}

\PROOF
	The relative entropy functional $ \RELENT{\mu}{ \mu_{ \BARV}} $ \VIZ{REfunc}  
with  $\gamma(x) =\mu(x) $ in \VIZ{recon}  becomes
\begin{eqnarray*}
 \RELENT{\mu}{ \mu_{ \BARV} } &=& \E_{\mu}\left[ \log\frac{\mu}{\mu_{ \BARV} } \right]=\int \mu(x) \log \frac{\mu(x)}{\mu_{ \BARV} (x)}  dx\\
&=&  \int\int \mu(x) \log \frac{\mu(x)}{\BARMV(\bx)\nu(x|\bx)}  d(x|\bx) d\bx =  \int\int \mu(x) \log \frac{\mu(x) \BARM(\bx)}{\BARMV(\bx) \mu(x)}  d(x|\bx) d\bx\\
&=& \int\int \mu(x) d(x|\bx)  \log \frac{\BARM(\bx)}{\BARMV(\bx) } d\bx
=\int \BARM(\bx) \log \frac{\BARM(\bx)}{\BARMV(\bx) } d\bx\\
&=& \E_{\BARM }\left[ \log \frac{\BARM(\bx)}{\BARMV (\bx)}\right]\COMMA
\end{eqnarray*}
thus, since $\BARM(\bx) =e^{-\beta \PMF(z)}/Z$ and $\BARMV (\bx) =  e^{-\beta   \BARV}/Z_{ \BARV}$,
\begin{equation*} 
 \RELENT{\mu}{\mu_{ \BARV} } = \E_{\BARM }\left[   -\beta \left(  \PMF  -   \BARV \right)\right] + \log Z_{ \BARV}/ Z  \PERIOD
\end{equation*}
 Expanding the logarithm and the exponential in the partition function term, when \\ $\PMF(\bx) - ~ \BARIT \IP(\bx)$ is small, we get
 \begin{eqnarray*}
 \log Z_{ \BARV}/ Z   &= & Z_{ \BARV}/ Z  -1 + \frac12 \left( Z_{ \BARV}/ Z  -1\right)^2 + \BIGO((Z_{ \BARV}/ Z -1)^2)  \COMMA \text{ and }\\
 Z_{ \BARV}/ Z  -1 &= &\E_{\BARM }\left[ e^{-\beta( \BARV - \PMF)}\right] -1 \\
&=& -\beta\E_{\BARM }\left[ ( \BARV - \PMF) \right]+\frac12\beta^2\E_{\BARM }\left[ ( \BARV - \PMF)^2 \right]  + \BIGO(\beta^3\E_{\BARM }\left[ ( \BARV - \PMF)^3\right] ) \PERIOD
 \end{eqnarray*}
Therefore
 \begin{eqnarray*}
 \log Z_{ \BARV}/ Z  =  -\beta\E_{\BARM }\left[ ( \BARV - \PMF) \right]+ \beta^2\E_{\BARM }\left[ ( \BARV - \PMF)^2 \right]  + \BIGO(\beta^3\E_{\BARM }\left[ ( \BARV - \PMF)^3\right]  ) \COMMA
 \end{eqnarray*}
 and the relative entropy is
 \begin{equation*} 
 \RELENT{\mu}{\mu_{ \BARV} } = \beta^2\E_{\BARM }\left[    \left(  \PMF  -   \BARV  \right)^2\right]  + \BIGO(\beta^3\E_{\BARM }\left[ | \BARV - \PMF|^3\right]  )  \PERIOD
\end{equation*}
We contrast this finding with the  minimization in $\Ll$ below.
For the  functional $\Ll(G;h)$ the representation \VIZ{chirep} holds
\begin{equation*} 
\Ll(G;h) = \Ll(F;h) + \beta^2 \E_{\mu}\left[ \left\| F(\COP)-G(\COP )\right\|^2 \right]\COMMA
\end{equation*}
thus    under the assumption  that   $ F(\bx) = \E_\mu\left[ h|\bx\right] = \FPMF(\bx)$, and that $G\in \Ee$ is of the form
$G_i(\bx) = -\nabla_{\bx_i} \BARIT \IP(\bx) $ for $  \BARV \in \Vv $  and  the definition of mean force $ \FPMF_i(\bx) = -\nabla_{\bx_i} \PMF(\bx),\ i=1,\dots,M  $ we have 
\begin{equation*} 
\Ll(G;h) = \Ll(\FPMF) + \beta^2 \E_{\mu}\left[ \left\|  \nabla \PMF( \COP)- \nabla   \BARIT \IP( \COP)\right\|^2 \right]\COMMA
\end{equation*}
thus, since $\E_{\mu}\left[ \phi( \xi)\right] = \E_{\BARM}\left[ \phi\right]   $ for any observable $\phi$ in $\R^{3M}$ we have 
\begin{equation*} 
\Ll(G;h) = \Ll(\FPMF) +  \beta^2\E_{\BARM}\left[ \left\|  \nabla \PMF(\bx) - \nabla   \BARIT \IP(\bx)\right\|^2 \right]\PERIOD
\end{equation*}
\hfill $\square$

Observing relations \VIZ{RERphi} and \VIZ{chiphi} we notice that the leading term at
the relative entropy approach   minimizes the average of the square of  potential difference
$ \left( \PMF  -   \BARV \right)^2$, i.e., it is an $L^2(\BARM)$ error. On the other hand, the   force matching minimizes the  average  of 
$\|  \nabla  \left(  \PMF  -   \BARV  \right) \|^2$,  an $H^1(\BARM)$ error, where
 $  H^1(\BARM) \!= \!\{g\in L^2(\BARM): \! \text{ weak first derivatives \! } Dg\in  L^2(\BARM) \}$.
 Thus, assuming the minimization problems have unique optimal solutions,  $ \BARV_{RE}^*$ and $ \BARV_{FM}^* $ for the relative entropy and force matching methods respectively, these solutions  differ by a constant.

\subsection{Structural based  parametrization methods.}\label{IBIM}
This section concerns an alternative family of   CG effective potentials  given by the  structure based or correlation based methods such as  the inverse Boltzmann, direct\cite{tsop1} and   iterative\cite{Soper1996}, and the inverse Monte Carlo  methods \cite{LyubLaa2004}.  
Theoretically, if one can   compute the  $n$-body  correlation function  $\BARIT{g}^{n}(z)$, $n<M$  from  the microscopic system simulations then according to the  relation\cite{mcquarrie2000statistical}
\begin{equation*}  
 \PMF(z^{(n)}) =  -\frac{1}{\beta} \log \BARIT{g}^{n}(z^{(n)}) 
 \end{equation*}
 where $z^{(n)}=(z_1,\dots,z_n) $,  $\BARIT{g}^{(n)}(z^{(n)}) $ is the $n$-body  correlation function, 
 the computation of $\PMF$ is straightforward, and the structural based methods, in principle, can provide exactly  the potential of the mean force, as is the case of the relative entropy and force matching methods.

  However,   the computation of  $\BARIT{g}^{n}(z) $  is not feasible for large $n$, and what is in practice  used at inverse Boltzmann and inverse Monte Carlo methods is the pair correlation 
\begin{equation*}
  \BARIT{g}^{(2)}(z_1,z_2) = \frac{(M-1)M}{ \rho^2}\int_{\R^{M-2}} \BARM(z) dz_{3}\dots dz_M\PERIOD
 \end{equation*}
 In homogeneous systems $\BARIT g^{(2)}(z_1,z_2)$ depends on the relative position between two particles $r=\|z_1-z_2\|$, $\BARIT g^{(2)}(r)$,  called the radial distribution function  
 $$ \BARIT{g}^{(2)}(r) =   \frac{(M-1)M}{ \rho^2} E_\mu[ \mathbbm{1}_{B(\bx_2,r)}(\bx_1)  |\bx] = \frac{(M-1)M}{ \rho^2}  \int_{\{x:\COP(x)=\bx\} } \mathbbm{1}_{B(\bx_2,r)}(\bx_1) \mu(x) dx\COMMA $$
that is the average density of finding the CG particle 1 at a distance $r$ from the particle  
2.  
Moreover, all structure based   methods rely  on Henderson's uniqueness theorem  \cite{Henderson1974}, which states that  for a given radial distribution function    there is a unique,   up to a constant, pair potential   $v(r)$ such that
$$  \PMF(z) = \sum_{i,j} v(\|z_i-z_j\|)\PERIOD$$

The structure based methods with the use of the pair radial distribution function in principle are comparable to the force matching and relative entropy when  the later ones consider the family of proposed   potentials, $\Vv$ in Theorem~\ref{comREFM}, to consist of pair interaction potentials.  The numerical comparison of all methods for molecular systems under equilibrium and non-equilibrium conditions is the subject of the future work\cite{HKKP2}.

\section{Discussion and conclusions}\label{Conclude}
  The main goal of all systematic CG approaches, based on statistical mechanics, is in principle to derive effective CG interactions as a numerical approximation of the many-body potential of the mean force, which for realistic molecular complex systems cannot be calculated exactly.
  
In this work we have presented a general   formalism for the development of CG methodologies for molecular systems.  
 Below we summarize the main outcomes of the detailed analysis presented in the previous sections: 

(a) The probabilistic formalism discussed   allows us to define a systematic   force matching, as a  CG minimization problem  both for linear and nonlinear CG maps.  This probabilistic formulation gives  a geometric representation of the force matching method, as is schematically depicted in Figure~\ref{scheme2}.
(b) A practical outcome of (a) is  the   connection of force matching with thermodynamic integration that provides a way on how to construct a 
 local mean force in order to best approximate the potential of mean force with force matching. Specifically, this connection introduces a family of  corresponding (to the CG map) coarsening transformations of the  microscopic forces (local mean force) that ensure  the best approximation of the PMF. This approach extends the work in ref. \cite{Voth2008a, Noid2011}, for any nonlinear CG map.
(c) CG methods based on relative entropy and force matching are in principle asymptotically equivalent, in the sense of Theorem~\ref{comREFM}, both for the case of linear and nonlinear coarse-graining maps. 
Furthermore  we prove, for  linear CG maps in a specific example of a system with $N$ molecules,   that the (un-weighted) total force exerted at each CG particle satisfies the force matching condition  when each particle is contributing to a  single CG particle, see the example~\ref{total}. This fact, along with the example of the nonlinear CG map studied in  Section~\ref{ReactionC}, suggest that for  complicated  linear and nonlinear  CG mappings   one can use appropriately  formula \VIZ{genhPMFW} and achieve the best approximation of PFM with the force matching method. 

Current work concerns the extension of this formalism, following the results in  ref., \cite{KP2013}, to coarse graining in  non-equilibrium systems, an important challenge where in principle CG methods fail.\cite{Baig2010}
The numerical  implementation of the formalism to different complex molecular systems\cite{HKKP2} is also the subject of  current studies.

\appendix
\section{Conditional expectation and coarse graining}\label{sigma_algebra}\label{AppendixCond}
 Let $(\R^{3N}, \Gg,\mu)$  be the probability space  induced  by the  random variable $X$ of atomic configuration.   $\Gg $ is the $\sigma$-algebra generated by the random variable $X$, i.e. it is the collection 
 $ \Gg  = \{ A\in \R^{3N} : \exists B \text{ Borel in } \R^{3N} \text{ s.t. } X^{-1}(B) = A\}\COMMA$
  Consider the coarse-grained random variable $\COP= \COP(X)$ and 
define the sub $\sigma$-algebra of $\Gg$, induced by $\COP$,
$$ \Gg_\COP = \{ A\in \R^{3N} : \exists C \text{ Borel in } \R^{3M} \text{ s.t. } \COP^{-1}(C) = A\}\COMMA$$
i.e.  any function $\phi: \R^{3N}\to \R$ that is $\Gg_\COP$-measurable is of the form 
$$\phi (x) = \phi(\COP (x))\PERIOD $$
Denote $\Omega(\bx) = \{ x\in \R^{3N}: \COP (x)= \bx\}$, the sub-manifold of $\R^{3N}$ corresponding to configurations $x$ at a fixed value of the coarse grained variable $\bx\in \R^{3M}$.
The conditional expectation with respect to $\Gg_\COP$ is the random variable $\E_\mu\left[\phi|\COP\right]$, defined by  
\begin{eqnarray*}
\E_\mu\left[\phi|\COP=z\right]=\E_\mu\left[\phi|\Gg_\COP\right](\bx)  = \frac{1}{\BARM(\bx)}\int_{\Omega(\bx)} \phi(x)  \mu(x)dx,\ \text{ for any } \bx\COMMA
\end{eqnarray*}
and for any $\Gg$-measurable $\phi$, with 
$$\BARM(\bx) = \int_{\Omega(\bx)}\mu(x)dx\COMMA$$ that is the average of $\phi$ keeping $\bx$ fixed. 

\section{Proofs}\label{appendix1}
\subsection{Proof of Lemma~\ref{lemmachi} }
Let  $F(\bx) =\E_{\mu}[h|\bx] $ then for  any $G\in \Ee$ holds
\begin{eqnarray*}
\Ll(G;h) &=& \E_\mu\left[ \left\| h  - G( \COP) \right\|^2 \right]  =\E_{\mu}\left[ \left\| h  - F(\COP )+F(\COP) -G( \COP) \right\|^2 \right]\\
&=& \E_{\mu}\left[ \left\| h  - F(\COP ) \right\|^2 \right]+\E_{\mu}\left[ \left\| F( \COP)-G(\COP) \right\|^2 \right] + 2 \E_{\mu}\left[ (h - F(\COP))(F(\COP)-G( \COP)) \right]\\
&=&  \E_{\mu}\left[ \left\|h  - F(\COP)\right\|^2 \right]+\E_{\mu}\left[ \left\| F(\COP)-G(\COP) \right\|^2 \right] \\
&=& \Ll(F;h) + \E_{\mu}\left[ \left\| F(\COP)-G(\COP )\right\|^2 \right] \COMMA
\end{eqnarray*}
since 
\begin{eqnarray*}
\E_{\mu}\left[ (h  - F(\COP))(F(\COP)-G(\COP)) \right] &=&  \E_{\mu}\left[ \E_{\mu}\left[ (h  - F( \COP))(F(\COP)-G(\COP)) \right] | \bx\right]\\
&=&  \E_{\mu}\left[ (F(\COP)-G( \COP)) \E_{\mu}\left[ (h  - F( \COP)) |\bx\right] \right] \\
&=&  \E_{\mu}\left[ (F( \COP)-G( \COP))\left( \E_{\mu}\left[ h  |\bx\right] - F(\COP ) \right)\right] = 0 \PERIOD
\end{eqnarray*}
Thus 
\begin{equation*}
\inf_{G\in \Ee} \Ll(G;h) = \Ll(F;h)\COMMA
\end{equation*}
and 
\begin{equation*}
  \Ll(G;h) = \Ll(F;h) +  \E_{\mu}\left[ \left\| F(\COP)-G( \COP ) \right\|^2 \right] \PERIOD
\end{equation*}
\hfill $\square$

\subsection{Proof of Theorem \ref{ghPMF} }
Let  the sub-manifold  $\Omega(\bx) = \{ x\in \R^{3N}: \COP (x)= \bx\}$ of $\R^{3N}$,   have the co-dimension $3M$, i.e., $\dim(\R^{3N})-\dim(\Omega(\bx)) = 3M$. The $\delta$   measure is defined as follows, for   any smooth test function $\phi:\R^{3N}\to \R$,
  \begin{equation*}
  \int_{\R^{3N}} \delta(\COP (x) - \bx ) \phi(x) dx =
  \int_{\Omega(\bx)} \phi \left(\text{det} \BFJ{\COP}\right)^{-1/2}d\Sigma_{\Omega(\bx)}\COMMA
  \end{equation*} 
 where    $\cdot^t$   denotes the matrix transpose, $\text{det}(\cdot)$  the matrix  determinant  and $\Sigma_{\Omega(\bx)}$ denotes the surface measure on $\Omega(\bx)$.   
  Let the mollifier on $\R^{3N}$
  $$
  \delta_{\epsilon}( \COP (x) - \bx) = \frac{1}{(2\pi \epsilon)^{3M/2}}e^{-\frac{1}{2\epsilon} |\COP (x) - \bx|^2}  \text{ for any  }\epsilon>0\PERIOD
  $$
 We have 
 \begin{eqnarray*}
 \nabla_{x_j} \delta_{\epsilon} ( \COP (x) - \bx)  &=& \sum_{i=1}^M  \nabla_{i} \delta_{\epsilon} ( \COP (x) - \bx)  \nabla_{x_j} \COP_i(x)\\
 &=& - \sum_{i=1}^M \nabla_{\bx_i} \delta_{\epsilon} ( \COP (x) - \bx)  \nabla_{x_j} \COP_i(x)\COMMA
 \end{eqnarray*}
 recalling the notation $(\BFD\COP)_{ij} = \nabla_{x_j} \COP_i $,  $i=1,\dots,M, j=1,\dots,N$, then we can write
 $$\nabla_x \delta_{\epsilon}( \COP (x) - \bx)  = - \BFD\COP^t \nabla_{\bx} \delta_{\epsilon}( \COP (x) - \bx) \COMMA $$
 and, since we assume that  $\BFW \BFD\COP^t(x) $ is invertible, we can write
 \begin{equation*}
 \BFW\nabla_x \delta_{\epsilon} ( \COP (x) - \bx) = -\BFW \BFD\COP^t \nabla_{\bx} \delta_{\epsilon}( \COP (x) - \bx) 
 \end{equation*}
 from which we have 
  \begin{equation*}
 \nabla_{\bx} \delta_{\epsilon} ( \COP (x) - \bx) = - \left(\BFW \BFD\COP^t\right)^{-1} \BFW \nabla_x \delta_{\epsilon} ( \COP (x) - \bx) \PERIOD
 \end{equation*}
 Taking the limit as $\epsilon \to 0$, in view of Lemma~\ref{lemmadeltae} in Appendix \ref{appendix1}, we have
 \begin{equation}\label{deltaderivative}
 \nabla_{\bx} \delta(\COP (x) -\bx) = - \left(\BFW \BFD\COP^t\right)^{-1} \BFW \nabla_x \delta(\COP (x) -\bx)\PERIOD
 \end{equation}
 We recall the definition of the potential of mean force \VIZ{PMF},
 $$
 \PMF(\bx) =  -\frac{1}{\beta} \log \BARM(\bx)  -\frac{1}{\beta} \log Z=   -\frac{1}{\beta} \log \int_{\Omega(\bx)} e^{-\beta \IP(x)} dx  \COMMA
 $$
 which we rewrite as
 \begin{equation*}
  \PMF(\bx) = -\frac{1}{\beta} \log \int_{\R^{3N}} \delta(\COP (x)-\bx)e^{-\beta \IP(x)}  dx  \PERIOD
 \end{equation*}
 Therefore, in view of the relation~\VIZ{deltaderivative}
{\small  \begin{eqnarray*}
 &&\nabla_{\bx} \PMF(\bx) =  -\frac{1}{\beta}   \frac{1}{\BARM(\bx) } \nabla_{\bx}   \int_{\R^{3N}} \delta(\COP (x)-\bx) e^{-\beta \IP(x)} dx\\
 &=&  -\frac{1}{\beta}   \frac{1}{\BARM(\bx) }  \int_{\R^{3N}}\nabla_{\bx}  \delta(\COP (x)-\bx)e^{-\beta \IP(x)} dx =  -\frac{1}{\beta}  \frac{1}{\BARM(\bx) } \lim_{\epsilon\to 0}    \int_{\R^{3N}}\nabla_{\bx}  \delta_{\epsilon}(\COP (x)-\bx)e^{-\beta \IP(x)} \\
 &=&   -\frac{1}{\beta}  \frac{1}{\BARM(\bx) }  \lim_{\epsilon\to 0}   \int_{\R^{3N}}\left[- \left(\BFW \BFD\COP^t\right)^{-1} \BFW \nabla_x \delta_{\epsilon}(\COP (x)-\bx)\right]e^{-\beta \IP(x)} dx \\
 &=&\!\!\! \frac{1}{\beta}  \! \frac{1}{\BARM(\bx) }    \! \lim_{\epsilon\to 0}    \!\!  \int_{\R^{3N}}\!\!\!\left[- \nabla_x \cdot[\left(\BFW \BFD\COP^t\right)^{-1}\BFW ] + \beta \left(\BFW \BFD\COP^t\right)^{-1} \BFW \nabla_x \IP(x)  \right]e^{-\beta \IP(x)}\!\delta_{\epsilon}(\COP (x)-\bx) dx \\
 &=&  \frac{1}{\beta}   \frac{1}{\BARM(\bx) }       \int_{\R^{3N}}\left[\beta \left(\BFW \BFD\COP^t\right)^{-1} \BFW \nabla_x \IP(x)  - \nabla_x \cdot[\left(\BFW \BFD\COP^t\right)^{-1} \BFW ]   \right]e^{-\beta \IP(x)}\delta (\COP (x)-\bx) dx\\
 &=&  \E_\mu\left[  \left(\BFW \BFD\COP^t\right)^{-1} \BFW \nabla_x \IP(x)  - \frac1\beta\nabla_x \cdot[\left(\BFW \BFD\COP^t\right)^{-1}\BFW ]  |\bx\right] \PERIOD
 \end{eqnarray*}
 }
Thus we conclude that 
\begin{equation*}
\FPMF(\bx) = -\nabla_{\bx} \PMF(\bx) =  \E_\mu\left[  \left(\BFW \BFD\COP^t\right)^{-1} \BFW \left(-\nabla_x \IP(x)  \right)+ \frac1\beta\nabla_x \cdot[\left(\BFW \BFD\COP^t\right)^{-1} \BFW ]  |\bx\right]\COMMA
\end{equation*}
that is 
\begin{equation*}
\FPMF(\bx) =   \E_\mu\left[  \left(\BFW \BFD\COP^t\right)^{-1} \BFW f(x)+ \frac1\beta\nabla_x \cdot[\left(\BFW \BFD\COP^t\right)^{-1} \BFW ]  |\bx\right]\PERIOD
\end{equation*}
  \hfill $\square$

  \begin{lem}\label{lemmadeltae}
  Let the mollifier on $\R^{3N}$ be
  $$
  \delta_{\epsilon}( \COP (x) - \bx) = \frac{1}{(2\pi \epsilon)^{3M/2}}e^{-\frac{1}{2\epsilon} |\COP (x) - \bx|^2}  \text{ for any  }\epsilon>0
  $$
   then for any test function $\phi:\R^{3N}\to \R$
   \begin{equation*} 
  \int_{\R^{3N}} \delta_{\epsilon}(\COP (x)-\bx) \phi(x) dx \to  \int\delta(\COP (x)-\bx) \phi(x) dx \ \ \text{ as } \epsilon \to 0\PERIOD
   \end{equation*}
   Furthermore,
   \begin{equation*}
  \int_{\R^{3N}} \nabla_x\delta_{\epsilon}(\COP (x)-\bx) \phi(x) dx \to  \int \nabla_x\delta(\COP (x)-\bx) \phi(x) dx \ \ \text{ as } \epsilon \to 0\PERIOD
   \end{equation*}
    
\end{lem}
 \PROOF  Let the smooth  test function $\phi:\R^{3N}\to \R$. We have 
 \begin{eqnarray*}
   \int_{\R^{3N}} \delta_{\epsilon}(\COP (x)-\bx) \phi(x) dx & = &\int_{\R^{3N}} \frac{1}{(2\pi \epsilon)^{3M/2}}e^{-\frac{1}{2\epsilon} |\COP (x) - \bx|^2}\phi(x) dx  \\
   &=& \frac{1}{(2\pi \epsilon)^{3M/2}}\int_{\R^{3M}}\int_{\Omega(\bx)} e^{-\frac{1}{\epsilon} |\COP (x) - \bx|^2}\phi(x) d\Sigma_{\Omega(\bx)} d\bx \PERIOD
 \end{eqnarray*}
Next we define    the  orthogonal projection onto $ \Omega(\bx)$, 
\begin{eqnarray*}
P_{\Omega(\bx)}: \R^{3N} \to \Omega(\bx), \quad
 x \mapsto  \xi \in \Omega(\bx), \ \xi = P_{\Omega(\bx)} x
 \end{eqnarray*}
 then 
$ \R^{3N} = \BFT_{\xi} \Omega(\bx) \oplus N_{\xi} \Omega(\bx) $, where $\BFT_{\xi}  $ and $N_{\xi} $ denote the tangent and normal space to $\Omega(\bx) $ at $\xi$. 
We denote the local coordinates on $\Omega(\bx)$ at $\xi$, $(\tau_1(\xi), \dots, \tau_{3N-3M}(\xi))$, then
$$x=x(\tau, \eta) = \xi(\tau) + \sum_{i=1}^{3M} \eta_i e_i(\xi(\tau))\COMMA$$
where $ \xi(\tau)\in \Omega(\bx)$ and $\sum_{i=1}^{3M} \eta_i e_i(\xi(\tau))$ is its normal conjugate,
therefore 
$$ dx = \left(\text{det} \BFJ(\tau)\right)^{ 1/2} d\tau d\eta$$
where   $\BFJ(\tau)= \nabla\xi(\tau)^t \nabla \xi(\tau)$ is the metric induced by the embedding of $\Omega(\bx)$. In other words 
$$ 
 \left(\text{det} \BFJ(\tau)\right)^{ 1/2} d\tau d\eta = d\Sigma_{\Omega(\bx)} \PERIOD
$$
Thus based on the expansion of $\frac12|\COP (x) - \bx|^2 $ around its minimum  on $ \Omega(\bx)$
$$ \frac12|\COP (x) - \bx|^2 = \frac12|\BFD\COP(x)(x - \bx)|^2 + \BIGO(|x-\bx|^4) \COMMA$$
 where $\BFD\COP(x)$ is the   matrix    with elements $(\BFD\COP)_{ij}(x) =  \partial_{x_j} \COP_i(x) , \ i=1,\dots,3M, j=1,\dots,3N$,   we can write 
\begin{eqnarray*}
&&\int_{\R^{3N}}\frac{1}{(2\pi \epsilon)^{3M/2}}e^{-\frac{1}{2\epsilon} |\COP (x) - \bx|^2}\phi(x) dx = \frac{1}{(2\pi \epsilon)^{3M/2}}\int_{\R^{3N}} e^{-\frac{1}{2\epsilon} |\BFD\COP(x)(x - \bx)|^2}\phi(x) dx + O(\epsilon)\\
&&=\frac{1}{(2\pi \epsilon)^{3M/2}}\int_{\R^{3N}} e^{-\frac{1}{2\epsilon} |\BFD\COP(x)\eta|^2}\left(\phi(\bx) +\phi'(\bx)\eta + O(|\eta|^2)  \right) \left(\text{det} \BFJ(\tau)\right)^{ 1/2} d\tau d\eta + O(\epsilon)\\
&& = \int_{\Omega(\bx)} \phi(\xi(\tau)) \left(\text{det} \BFJ{\COP}(\xi(\tau))\right)^{-1/2} \left(\text{det} \BFJ(\tau)\right)^{ 1/2} d\tau d\eta + O(\epsilon)\\
&&\to   \int_{\Omega(\bx)} \phi \ \ \left(\text{det} \BFJ{\COP} \right)^{-1/2} d\Sigma_{\Omega(\bx)} =  \int_{\R^{3N}} \delta(\COP (x) - \bx ) \phi(x) dx \ \  \text{ as } \epsilon \to 0 \PERIOD 
 \end{eqnarray*}
  Furthermore, 
   \begin{equation*}
  \int_{\R^{3N}} \nabla_x\delta_{\epsilon}(\COP (x)-\bx) \phi(x) dx \to  \int \nabla_x\delta(\COP (x)-\bx) \phi(x) dx \ \ \text{ as } \epsilon \to 0\COMMA
   \end{equation*}
  which is proved similarly.
 \hfill $\square$



\begin{acknowledgments}
The research of E.K. and V.H. was supported by the European Union (ESF) and Greek national funds through the Operational Program “Education and Lifelong Learning” of the NSRF-Research Funding Program: THALIS. The research of M.K. was supported in part by the Office of Advanced Scientific Computing Research, U.S. Department of Energy under Contract No. DE- SC0010723. The research of P.P. was partially supported by the U.S. Department of Energy Office of Science, Office of Advanced Scientific Computing Research, Applied Mathematics program under Award Number DE-SC-0007046. 
\end{acknowledgments}


%

\bibliography{reference1}

\end{document}